\documentclass{jcmlatex}
\def\JCMvol{xx}
\def\JCMno{x}
\def\JCMyear{200x}
\def\JCMreceived{xxx}
\def\JCMrevised{xxx}
\def\JCMaccepted{xxx}

\usepackage{bm}
\usepackage{color}
\usepackage[hang]{subfigure}
\usepackage{multirow}\usepackage{epstopdf}

\newcommand\bbR{\mathbb{R}}
\newcommand\dd{\mathrm{d}}

\newcommand\bx{\bm{x}}
\newcommand\intt{\int_{t^n}^{t^{n+1}}}
\newcommand\ph{{\varphi_h}}
\newcommand\il{{i-\frac12}}
\newcommand\ir{{i+\frac12}}
\newcommand\jl{{j-\frac12}}
\newcommand\jr{{j+\frac12}}

\newcommand\pd[2]{\dfrac{\partial {#1}}{\partial {#2}}}

\newcommand\abs[1]{\lvert #1 \rvert}

\newcommand\pro[2]{\left<#1,#2\right>}
\newcommand{\tabincell}[2]{\begin{tabular}{@{}#1@{}}#2\end{tabular}}

\begin{document}

\markboth{J.M. Duan AND H.Z. Tang}{ADER-DG scheme for Hamilton-Jacobi equations}

\title{An efficient ADER discontinuous Galerkin scheme for
directly solving Hamilton-Jacobi equation\footnote{It is an invited paper for the special issue
on ``The Sixth Chinese-German Workshop on Computational and Applied Mathematics''.}}

\author{Junming Duan
\thanks{LMAM, School of Mathematical Sciences, Peking University,
Beijing 100871,  China 
}
\and
Huazhong Tang
\thanks{HEDPS, CAPT \& LMAM, School of Mathematical Sciences, Peking University,
Beijing 100871; School of Mathematics and Computational Science,
 Xiangtan University, Hunan Province, Xiangtan 411105, P.R. China}}

\maketitle

\begin{abstract}
This paper  proposes an efficient ADER (Arbitrary DERivatives in space and time) discontinuous Galerkin (DG) scheme to directly solve the Hamilton-Jacobi equation.
Unlike multi-stage Runge-Kutta methods used in the Runge-Kutta DG (RKDG) schemes,
the ADER scheme is one-stage in time discretization, which is desirable in many applications.
The ADER scheme used here relies on a local continuous spacetime Galerkin predictor
instead of the usual Cauchy-Kovalewski procedure to achieve high order accuracy both in space and time.
In such predictor step,  a local Cauchy problem in each cell
is solved based on a weak formulation of the original equations in spacetime.
The resulting spacetime representation of the numerical solution provides the
temporal accuracy that matches the spatial accuracy of the underlying DG solution.
The scheme is formulated in the modal space and the volume integral and the numerical fluxes at the cell interfaces can be explicitly written.
The explicit formulae of the scheme at third order is provided on two-dimensional structured meshes.
The computational complexity of the ADER-DG scheme is compared to that of  the RKDG scheme.
Numerical experiments are also provided to demonstrate the accuracy and efficiency of our scheme.
\end{abstract}



\section{Introduction}
Consider the Hamilton-Jacobi (HJ) equation
\begin{equation}
  \varphi_t+H(\nabla_{\bx}\varphi,\bx)=0, \quad
  \varphi(\bx,0)=\varphi^0(\bx), \quad \bx\in\Omega\in\bbR^d,
  \label{eq:HJ}
\end{equation}
with suitable boundary conditions, where $H(\cdot)$ denotes the Hamiltonian.
The HJ equations are used in many application areas, such as optimal control theory, geometrical optics,
crystal growth, image processing and computer vision.
The solutions of such equations are continuous
but their derivatives could be discontinuous even if the initial condition is smooth.
Viscosity solutions were firstly introduced and studied
in \cite{Crandall1983,Crandall1984},
which are the unique physically relevant solutions.

It is well known that the HJ equations are closely related to hyperbolic conservation laws,
thus many successful numerical methods for the conservation laws can be adapted for solving the HJ equations.
In \cite{Crandall1984}, a monotone finite difference scheme was introduced and proved to be convergent to the viscosity solution.
A second order finite difference essentially non-oscillatory (ENO) scheme was developed in \cite{Osher1988},
and then a higher-order weighted ENO (WENO) scheme is proposed in \cite{Jiang2000}.
Tang and his collaborators \cite{Tang2003} developed an adaptive mesh redistribution method for the HJ equations in two- and three-dimensions.
Qiu et al. \cite{Qiu2005} developed the Hermite WENO (HWENO) schemes of the HJ equations.
The high order finite difference WENO scheme on unstructured meshes was developed in \cite{Zhang2003}, but its implementation is a bit complicated.

Alternatively,   a DG method was designed in \cite{HU1999} to solve the HJ equations,
and its reinterpretation and simplified implementation was given in \cite {Li2005}.
Those DG methods were based on the fact that the derivatives of the solution
satisfied the conservation laws. It was  correct  in the one-dimensional case
but at risk in the multi-dimensional case
because corresponding multi-dimensional conservation laws is  only weakly hyperbolic in general.
Later, a DG method for directly solving the HJ equations with convex Hamiltonians was proposed in \cite{Cheng2007}. It was further improved and  a new DG method was derived
for directly solving the general HJ equations with nonconvex Hamiltonians in \cite{Cheng2014}.
This paper will construct the scheme based on the RKDG scheme in \cite{Cheng2014}.
The RKDG method \cite{Cockburn2001} was originally designed to solve conservation laws,
which has the advantages of flexibility on complicated geometries and  a compact stencil, and is easy to obtain high order accuracy.

Most of the above methods use the multi-stage Runge-Kutta time discretization,
thus have the advantage of simplicity  but are time-consuming because
at each stage, the volume integration and the numerical fluxes at cell interfaces
have to be calculated and  the nonlinear limiters should be performed
to suppress the numerical oscillations.
Thus, in order to save the computational cost,
it is desirable to use an alternative to the multi-stage Runge-Kutta method.
One choice is the Lax-Wendroff type time discretization, which
 converts all (or partial,
when approximations with certain accuracy are expected) time derivatives
in a temporal Taylor expansion of the solution into spatial derivatives
by repeatedly using the underlying differential equation and its differentiated forms \cite{Guo2011}.
In \cite{Guo2011}, a local-structure-preserving DG method with Lax-Wendroff type time discretization was proposed for solving the HJ equations. It is shown that such method is relatively more efficient than the RKDG method in \cite{Li2005}.
But the Cauchy-Kowalewski procedure may become a little cumbersome
when we want to construct a high order scheme.
This paper will use the time discretization (named ADER)
proposed in \cite{Dumbser2008a,Dumbser2008}.
The ADER scheme has  been successfully applied to the (magneto) hydrodynamics and relativistic (magneto) hydrodynamics with stiff or non-stiff source terms
\cite{Balsara2013,Balsara2009,Dumbser2008a,Dumbser2008,Fambri2018}.
It is based on a local spacetime Galerkin predictor step, at which  a local Cauchy problem
is solved in each cell, based on a weak formulation of the original partial differential equations in spacetime.
Through the above procedure, the resulting spacetime representation of the numerical solution provides the temporal accuracy that
matches the spatial accuracy of the underlying DG solution.
The ADER scheme is a one-step one-stage time discretization,
which means that  the volume integration
and the numerical fluxes terms at cell interfaces
are  only  calculated once at each time step.
Our ADER-DG scheme is formulated in modal space.
Thanks to the spacetime representation of the numerical solution,
 we can write down explicit formulae of the scheme
using the strategy presented   in \cite{Balsara2013},
and we will provide the implementation details of the scheme at third order
on two-dimensional structured meshes.
Our ADER-DG scheme can capture the viscosity solution accurately and efficiently, and
 will be validated by the analysis of the computational complexity and the numerical experiments.

The paper is organized as follows.
Section \ref{sec:ADER-DG}   presents the general formulation of our one- and two-dimensional ADER-DG schemes.
 Section \ref{sec:predictor}   introduces the local spacetime continuous Galerkin predictor,
and gives a detailed description of the two-dimensional predictor step at third order.
Section \ref{sec:flux} describes the calculation of the volume integration and the numerical fluxes terms at the cell interfaces,
and  the computational complexity of our ADER-DG scheme   will be compared to the RKDG scheme in \cite{Cheng2014}.
Section \ref{sec:numer} presents numerical experiments and the
concluding remarks are given in Section \ref{sec:conclusion}.

\section{General formulation of the ADER-DG schemes}\label{sec:ADER-DG}
This section  will present the general formulation of the ADER-DG schemes,
in which the numerical fluxes and the penalty terms adding to the numerical fluxes
are firstly developed in \cite{Cheng2014}.

\subsection{One-dimensional ADER-DG scheme}
Let us consider the one-dimensional HJ equation at first.
In this case, \eqref{eq:HJ} becomes
\begin{equation}
  \varphi_t+H(\varphi_x,x)=0, \quad \varphi(x,0)=\varphi^0(x).
  \label{eq:HJ1D}
\end{equation}
Assume the computational domain  $[a,b]$ is  divided into $N$ cells,
$I_i=(x_{i-\frac12},x_{i+\frac12})$, $i=1,2,\cdots,N$, where
\begin{equation*}
  a=x_{\frac12}<x_{\frac32}<\dots<x_{N+\frac12}=b.
\end{equation*}
Denote the  center of $I_i$ and the mesh size as
$x_i=\dfrac12(x_{i-\frac12}+x_{i+\frac12})$
and $\Delta x_i=x_{i+\frac12}-x_{i-\frac12}$ respectively.
Moreover, use $H_1=\pd{H}{\varphi_x}$ to denote the partial derivative of the Hamiltonian with respect to $\varphi_x$.

The spatial DG approximation space is
\begin{equation}
  V^k_h=\{v:v|I_i\in P^k(I_i),i=1,2,\dots,N\},
\end{equation}
where $P^k(I_i)$ denotes all polynomials of degree at most $k$ on $I_i$.
Assume the current time interval is $[t^n,t^{n+1}]$, and the time stepsize is $\Delta t=t^{n+1}-t^n$.
Following \cite{Cheng2014}, if
multiplying \eqref{eq:HJ1D} with the
test function $\phi_m(x)\in V^k_h,k\geqslant 1$,
introducing the numerical fluxes, adding penalty terms
for the numerical fluxes at the interfaces of computational cells,
and integrating it over the spacetime control volume $I_i\times [t^n,t^{n+1}]$,
then one has
\begin{align}
  &\intt\int_{I_i}\phi_m(x)(\partial_t\varphi(x,t)+H(\partial_x \varphi (x,t),x))\dd x\dd t \nonumber \\
  &+\intt\min(\tilde{H}_{1,\varphi}(x_\ir,t),0)[\varphi]_\ir(\phi_m)_\ir^-\dd t
  +\intt\max(\tilde{H}_{1,\varphi}(x_\il,t),0)[\varphi]_\il(\phi_m)_\il^+\dd t \nonumber \\
  &-\intt C\Delta x_i\left( S_{1,\varphi}(x_\ir,t)-\abs{\tilde{H}_{1,\varphi}(x_\ir,t)}\right)
  [(\varphi)_x]_\ir(\phi_m)_\ir^-\dd t \nonumber \\
  &-\intt C\Delta x_i\left( S_{1,\varphi}(x_\il,t)-\abs{\tilde{H}_{1,\varphi}(x_\il,t)}\right)
  [(\varphi)_x]_\il(\phi_m)_\il^+\dd t \nonumber \\
  &=0, \quad \forall i = 1,2,\dots,N,
  \label{eq:DG1D}
\end{align}
where $[u]=u^+-u^-$ denotes the jump of function $u$ at the cell interface,
  the superscripts $+,-$ denote the right, and left limits of a function,
$C$ is a positive penalty parameter chosen as $0.25$ in this paper,
and $\tilde{H}_{1,\varphi}$ and $S_{1,\varphi}$ are the Roe speed and
the parameter to identify the entropy violating cells \cite{Cheng2014}.
For all $t\in[t^n,t^{n+1}]$, assume $(x_*,t)$ is a point located at the cell interface,
then $\tilde{H}_{1,\varphi}$ and $S_{1,\varphi}$ are defined by
\begin{align}
  &\tilde{H}_{1,\varphi}(x_*,t)=\begin{cases}
    \dfrac{H(\varphi_x(x_*^+,t),x_*^+)-H(\varphi_x(x_*^-,t),x_*^-)}
    {\varphi_x(x_*^+,t)-\varphi_x(x_*^-,t)},
    \quad & \text{if} \quad \varphi_x(x_*^-,t)\not=\varphi_x(x_*^+,t),\\
    \dfrac12\left( H_1(\varphi_x(x_*^+,t),x_*^+)+H_1(\varphi_x(x_*^-,t),x_*^-)\right),
    \quad & \text{if} \quad \varphi_x(x_*^-,t)=\varphi_x(x_*^+,t),\\
  \end{cases}\nonumber \\
  &\delta_{1,\varphi}(x_*,t)=\max(0,\tilde{H}_{1,\varphi}(x_*,t)-H_1((\varphi)_x(x_*^-,t),x_*^-),
  H_1((\varphi)_x(x_*^+,t),x_*^+)-\tilde{H}_{1,\varphi}(x_*,t)),\nonumber \\
  &S_{1,\varphi}(x_*,t)=\max(\delta_{1,\varphi}(x_*,t),\abs{\tilde{H}_{1,\varphi}(x_*,t)}).
\end{align}
It is worth noting  that the above definitions only make
sense for $k\geqslant 1$.
That is why we choose the DG space as $V^k_h,k\geqslant 1$.

Calculate the time derivative parts in \eqref{eq:DG1D},
restrict the solutions $\varphi(x,t)$ to $V^k_h$,
thus \eqref{eq:DG1D} becomes
\begin{align}
  &\int_{I_i}\phi_m(x)(\ph(x,t^{n+1})-\ph(x,t^{n}))\dd x
  +\intt\int_{I_i}\phi_m(x)H(\partial_x q_h(x,t),x)\dd x\dd t \nonumber \\
  &+\intt\min(\tilde{H}_{1,q_h}(x_\ir,t),0)[q_h]_\ir(\phi_m)_\ir^-\dd t
  +\intt\max(\tilde{H}_{1,q_h}(x_\il,t),0)[q_h]_\il(\phi_m)_\il^+\dd t \nonumber \\
  &-\intt C\Delta x_i\left( S_{1,q_h}(x_\ir,t)-\abs{\tilde{H}_{1,q_h}(x_\ir,t)}\right)
  [(q_h)_x]_\ir(\phi_m)_\ir^-\dd t \nonumber \\
  &-\intt C\Delta x_i\left( S_{1,q_h}(x_\il,t)-\abs{\tilde{H}_{1,q_h}(x_\il,t)}\right)
  [(q_h)_x]_\il(\phi_m)_\il^+\dd t \nonumber \\
  &=0, \quad \forall i = 1,2,\cdots,N,
  \label{eq:ADER-DG1D}
\end{align}
where an element-local spacetime predictor solution $q_h(x,t)$ is introduced
to replace the solution $\varphi(x,t)$ in the integral of the Hamiltonian and the numerical fluxes,
which is a high order approximation polynomial obtained by using the
local spacetime Galerkin predictor,
and will be presented in detail in Section \ref{sec:predictor}.
At time $t^n$, the DG solution $\ph(x,t^n)$ is known, if we get $q_h(x,t)$,
then the DG solution can be evolved to time $t^{n+1}$ as $\ph(x,t^{n+1})$
from \eqref{eq:ADER-DG1D}.

\subsection{Two-dimensional ADER-DG scheme on structured meshes}
Consider the two-dimensional HJ equation
\begin{equation}
  \varphi_t+H(\varphi_x,\varphi_y,x,y)=0, \quad \varphi(x,y,0)=\varphi^0(x,y),
  \label{eq:HJ2D}
\end{equation}
and the computational domain  $[a,b]\times[c,d]$ is divided into $N_x\times N_y$ cells.
$I_{i,j}=J_i\times K_j$ with $J_i = [x_\il,x_\ir]$, $K_j=[y_\jl,y_\jr]$,
$\Delta x_i=x_\ir-x_\il$, and $\Delta y_j=y_\jr-y_\jl$,
$i=1,\dots,N_x$, $j=1,\dots,N_y$.
Use $H_1=\pd{H}{\varphi_x}$ and $H_2=\pd{H}{\varphi_y}$ to denote the partial derivatives of the Hamiltonian
with respect to $\varphi_x$ and $\varphi_y$ respectively.

The spatial DG approximation space is
\begin{equation}
  V^k_h=\{v:v|I_{i,j}\in P^k(I_{i,j}),i=1,2,\dots,N_x,j=1,2,\dots,N_y\},
\end{equation}
where $P^k(I_{i,j})$ denotes all polynomials of degree at most $k$ on $I_{i,j}$,~$k\geqslant 1$.
Assume the current time interval is $[t^n,t^{n+1}]$, and the time step is $\Delta t=t^{n+1}-t^n$.
Multiplying \eqref{eq:HJ2D} with test functions $\phi_m(x,y)\in V^k_h,k\geqslant 1$,
introducing the numerical fluxes and adding penalty terms for the numerical fluxes
at the cell interfaces, and then integrating it
over the spacetime control volume $I_{i,j}\times [t^n,t^{n+1}]$
can give
\begin{align}
  &\intt\int_{I_{i,j}}\phi_m(x,y)(\partial_t\varphi(x,y,t)+H(\partial_x \varphi
  (x,y,t),\partial_y \varphi (x,y,t),x,y))\dd x\dd y\dd t \nonumber \\
  &+\intt\int_{K_j}\min(\tilde{H}_{1,\varphi}(x_\ir,y,t),0)[\varphi](x_\ir,y,t)\phi_m(x_\ir^-,y)\dd y\dd t \nonumber\\
  &+\intt\int_{K_j}\max(\tilde{H}_{1,\varphi}(x_\il,y,t),0)[\varphi](x_\il,y,t)\phi_m(x_\il^+,y)\dd y\dd t \nonumber\\
  &+\intt\int_{J_i}\min(\tilde{H}_{2,\varphi}(x,y_\jr,t),0)[\varphi](x,y_\jr,t)\phi_m(x,y_\jr^-)\dd x\dd t \nonumber\\
  &+\intt\int_{J_i}\max(\tilde{H}_{2,\varphi}(x,y_\jl,t),0)[\varphi](x,y_\jl,t)\phi_m(x,y_\jl^+)\dd x\dd t \nonumber\\
  &-C\Delta x_i\intt\int_{K_j}\left(
  S_{1,\varphi}(x_\ir,y,t)-\abs{\tilde{H}_{1,\varphi}(x_\ir,y,t)}\right)
  [\varphi_x](x_\ir,y,t)\phi_m(x_\ir^-,y)\dd y\dd t \nonumber \\
  &-C\Delta x_i\intt\int_{K_j}\left(
  S_{1,\varphi}(x_\il,y,t)-\abs{\tilde{H}_{1,\varphi}(x_\il,y,t)}\right)
  [\varphi_x](x_\il,y,t)\phi_m(x_\il^+,y)\dd y\dd t \nonumber \\
  &-C\Delta y_j\intt\int_{J_i}\left(
  S_{2,\varphi}(x,y_\jr,t)-\abs{\tilde{H}_{2,\varphi}(x,y_\jr,t)}\right)
  [\varphi_y](x,y_\jr,t)\phi_m(x,y_\jr^-)\dd x\dd t \nonumber \\
  &-C\Delta y_j\intt\int_{J_i}\left(
  S_{2,\varphi}(x,y_\jl,t)-\abs{\tilde{H}_{2,\varphi}(x,y_\jl,t)}\right)
  [\varphi_y](x,y_\jl,t)\phi_m(x,y_\jl^+)\dd x\dd t \nonumber \\
  &=0, \quad \forall i = 1,2,\dots,N_x, j = 1,2,\dots,N_y.
  \label{eq:DG2D}
\end{align}
For all $t\in[t^n,t^{n+1}],\forall y$, if assuming $(x_*,y,t)$ is a point located at the cell interface in the $x$-direction,
then the Roe speed and the parameters to identify the entropy violating cells in the scheme
are given by
\begin{align*}
  &\tilde{H}_{1,\varphi}(x_*,y,t)\\
  &=\begin{cases}
    \dfrac{H(\varphi_x(x_*^+,y,t),\overline{\varphi_y},x_*^+,y)-H(\varphi_x(x_*^-,y,t),\overline{\varphi_y},x_*^-,y)}
    {\varphi_x(x_*^+,y,t)-\varphi_x(x_*^-,y,t)},
    &\varphi_x(x_*^-,y,t)\not=\varphi_x(x_*^+,y,t),\\
    \dfrac12\left(
    H_1(\varphi_x(x_*^+,y,t),\overline{\varphi_y},x_*^+,y)+H_1(\varphi_x(x_*^-,y,t),\overline{\varphi_y},x_*^-,y)\right),
    &\varphi_x(x_*^-,y,t)=\varphi_x(x_*^+,y,t),\\
  \end{cases}\\
  &\delta_{1,\varphi}(x_*,y,t)\\
  &=\max(0,\tilde{H}_{1,\varphi}(x_*,y,t)-H_1(\varphi_x(x_*^-,y,t),\overline{\varphi_y},x_*^-,y),
  H_1(\varphi_x(x_*^+,y,t),\overline{\varphi_y},x_*^+,y)-\tilde
{H}_{1,\varphi}(x_*,y,t)),\\
&S_{1,\varphi}(x_*,y,t)=\max(\delta_{1,\varphi}(x_*,y,t),\abs{\tilde{H}_{1,\varphi}(x_*,y,t)}),
\end{align*}
where $\overline{\varphi_y}=\dfrac12(\varphi_y(x_*^+,y,t)+\varphi_y(x_*^-,y,t))$ is the average of
the tangential derivative.

Similarly, for all $t\in[t^n,t^{n+1}],\forall x$,
 if denoting $(x,y_*,t)$ as a point located at the cell interface in the $y$-direction,
then the Roe speed and the parameters are given by
\begin{align*}
  &\tilde{H}_{2,\varphi}(x,y_*,t)\\
  &=\begin{cases}
    \dfrac{H(\overline{\varphi_x},\varphi_y(x,y_*^+,t),x,y_*^+)-H(\overline{\varphi_x},\varphi_y(x,y_*^-,t),x,y_*^-)}
    {\varphi_y(x,y_*^+,t)-\varphi_y(x,y_*^-,t)},
    &\varphi_y(x,y_*^-,t)\not=\varphi_y(x,y_*^+,t),\\
    \dfrac12\left(
    H_2(\overline{\varphi_x},\varphi_y(x,y_*^+,t),x,y_*^+)+H_2(\overline{\varphi_x},\varphi_y(x,y_*^-,t),x,y_*^-)\right),
    &\varphi_y(x,y_*^-,t)=\varphi_y(x,y_*^+,t),\\
  \end{cases}\\
  &\delta_{2,\varphi}(x,y_*,t)\\
  &=\max(0,\tilde{H}_{2,\varphi}(x,y_*,t)-H_2(\overline{\varphi_x},\varphi_y(x,y_*^-,t),x,y_*^-),
  H_2(\overline{\varphi_x},\varphi_y(x,y_*^+,t),x,y_*^+)-\tilde
{H}_{2,\varphi}(x,y_*,t)),\\
&S_{2,\varphi}(x,y_*,t)=\max(\delta_{2,\varphi}(x,y_*,t),\abs{\tilde{H}_{2,\varphi}(x,y_*,t)}),
\end{align*}
where $\overline{\varphi_x}=\dfrac12(\varphi_x(x,y_*^+,t)+\varphi_x(x,y_*^-,t))$.

After calculating the time derivative parts in \eqref{eq:DG2D}
and restricting the solutions $\varphi(x,y,t)$ to $V^k_h$,
then \eqref{eq:DG2D} becomes
\begin{align}
  &\int_{I_{i,j}}\phi_m(x,y)(\ph(x,y,t^{n+1})-\ph(x,y,t^{n}))\dd x\dd y \nonumber \\
  &+\intt\int_{I_{i,j}} \phi_m(x,y)H(\partial_x q_h(x,y,t),\partial_y q_h(x,y,t),x,y)\dd x\dd y\dd t \nonumber \\
  &+\intt\int_{K_j}\min(\tilde{H}_{1,{q_h}}(x_\ir,y,t),0)[{q_h}](x_\ir,y,t)\phi_m(x_\ir^-,y)\dd y\dd t \nonumber\\
  &+\intt\int_{K_j}\max(\tilde{H}_{1,{q_h}}(x_\il,y,t),0)[{q_h}](x_\il,y,t)\phi_m(x_\il^+,y)\dd y\dd t \nonumber\\
  &+\intt\int_{J_i}\min(\tilde{H}_{2,{q_h}}(x,y_\jr,t),0)[{q_h}](x,y_\jr,t)\phi_m(x,y_\jr^-)\dd x\dd t \nonumber\\
  &+\intt\int_{J_i}\max(\tilde{H}_{2,{q_h}}(x,y_\jl,t),0)[{q_h}](x,y_\jl,t)\phi_m(x,y_\jl^+)\dd x\dd t \nonumber\\
  &-C\Delta x_i\intt\int_{K_j}\left(
  S_{1,{q_h}}(x_\ir,y,t)-\abs{\tilde{H}_{1,{q_h}}(x_\ir,y,t)}\right)
  [(q_h)_x](x_\ir,y,t)\phi_m(x_\ir^-,y)\dd y\dd t \nonumber \\
  &-C\Delta x_i\intt\int_{K_j}\left(
  S_{1,{q_h}}(x_\il,y,t)-\abs{\tilde{H}_{1,{q_h}}(x_\il,y,t)}\right)
  [(q_h)_x](x_\il,y,t)\phi_m(x_\il^+,y)\dd y\dd t \nonumber \\
  &-C\Delta y_j\intt\int_{J_i}\left(
  S_{2,{q_h}}(x,y_\jr,t)-\abs{\tilde{H}_{2,{q_h}}(x,y_\jr,t)}\right)
  [(q_h)_y](x,y_\jr,t)\phi_m(x,y_\jr^-)\dd x\dd t \nonumber \\
  &-C\Delta y_j\intt\int_{J_i}\left(
  S_{2,{q_h}}(x,y_\jl,t)-\abs{\tilde{H}_{2,{q_h}}(x,y_\jl,t)}\right)
  [(q_h)_y](x,y_\jl,t)\phi_m(x,y_\jl^+)\dd x\dd t \nonumber \\
  &=0, \quad \forall i = 1,2,\dots,N_x, j = 1,2,\dots,N_y,
  \label{eq:ADER-DG2D}
\end{align}
where the element-local spacetime predictor solution $q_h(x,y,t)$ will be introduced in Section \ref{sec:predictor}.
The remaining task is to give $q_h(x,y,t)$.

\section{Local spacetime continuous Galerkin predictor}\label{sec:predictor}
Unlike the classical ADER schemes in  \cite{Titarev2002,Titarev2005}
using Cauchy-Kovalewski procedure, which may become cumbersome for high order schemes,
the new formulation of ADER schemes proposed in   \cite{Dumbser2008a}
is based on a local weak formulation of the governing PDE in spacetime.
The new ADER schemes rely on an iterative predictor step to obtain the spacetime representation
of the solution within each cell, i.e., the previous mentioned
 local spacetime predictor solution $q_h$.
This part will  construct the predictor step,
and give the implementation details
of the predictor step in the two-dimensional case at third order.

\subsection{General formulation of continuous Galerkin predictor}
For the sake of convenience, we will only consider the two-dimensional case.
Assume the spatial coordinates in the reference element is $(\xi,\eta)\in[-\frac12,\frac12]^2$,
and the temporal coordinates in the reference element is $\tau\in[0,1]$.
In the reference element, Eq. \eqref{eq:HJ2D} can be written as
\begin{equation}
  \pd{\varphi}{\tau}+h\left(\dfrac{1}{\Delta
  x}\pd{\varphi}{\xi},\dfrac{1}{\Delta y}\pd{\varphi}{\eta},\xi,\eta\right)=0,
  \label{eq:HJref}
\end{equation}
where $h=\Delta tH$, and $\Delta x,\Delta y$ are the mesh sizes of the cell.
The ADER scheme   used here is a modal variant of the ADER scheme
with a continuous Galerkin representation in time described in \cite{Dumbser2008a}.
Assume that there are $L$ spacetime basis functions in the reference element,
$\theta_l=\theta_l(\xi,\eta,\tau),l=0,\cdots,L-1$.
The continuous Galerkin approach  requires that the first $L_s$ elements in the set of basis functions
only depend on the space but not on time $\tau$,
that is to say, $\theta_l(\xi,\eta,\tau)$ only depend on the space, $l=0,\cdots,L_s-1$.
Now the numerical solution $q_h$  can be represented in the basis space as
\begin{equation}
  q_h(\xi,\eta,\tau)=\sum_{l=0}^{L-1}\hat{q_l}\theta_l(\xi,\eta,\tau),
  \label{eq:q_expansion}
\end{equation}
where
$\hat{q}\equiv(\hat{q}_0,\cdots,\hat{q}_{L_s-1},\hat{q}_{L_s},\cdots,\hat{q}_{L-1})^\mathrm{T}$
is a vector of modes.
Similarly, the Hamiltonian can also be represented in the form of \eqref{eq:q_expansion},
$\hat{h}\equiv(\hat{h}_0,\cdots,\hat{h}_{L_s-1},\hat{h}_{L_s},\cdots,\hat{h}_{L-1})^\mathrm{T}$.
The transcription from $\hat{q}$ to $\hat{h}$ will be given in the next subsection.
Another simplification of the continuous Galerkin approach is that
the solution $q_h(\xi,\eta,\tau)$ is continuous with the initial condition $\varphi_h^n(\xi,\eta)$ at $\tau=0$,
which means we only have to calculate $\hat{h}_l,l=0,\cdots,L_s-1$ once at $\tau=0$.
If the initial condition can be represented in the modal space as
\begin{equation}
  \varphi_h^n(\xi,\eta)=\sum_{l=0}^{L_s-1}\hat{w}_l\theta_l(\xi,\eta,\tau=0),
\end{equation}
then at $\tau=0$, $\hat{q}_l=\hat{w}_l,l=0,\cdots,L_s-1$.

Applying the Galerkin approach to \eqref{eq:HJref} gives
\begin{equation}
  \pro{\theta_k}{\pd{\theta_l}{\tau}}\hat{q}_l + \pro{\theta_k}{\theta_l}\hat{h}_l=0,
  \label{eq:CG}
\end{equation}
where the angled brackets denote the spacetime integration over the reference element,
and the Einstein summation convection is used.
Eq. \eqref{eq:CG} can be rewritten in the matrix-vector form
\begin{equation}
  K_\tau\hat{q} + M\hat{h}=0,
  \label{eq:CG_predictor}
\end{equation}
where $K_\tau$ and $M$ are the time-stiffness matrix and the mass matrix respectively,
and the $(k,l)$-th elements of them are
\begin{equation}
  K_{\tau;k,l}= \pro{\theta_k}{\pd{\theta_l}{\tau}},
  M_{k,l}=\pro{\theta_k}{\theta_l}.
\end{equation}
From the previous assumption and simplification, we know that only the last $L-L_s$ elements
of $\hat{q}$ are needed to be determined in the continuous Galerkin predictor step.
So we can split $\hat{q}$ into two parts $\hat{q}=(\hat{q^0},\hat{q^1})^\mathrm{T}$,
where $\hat{q}^0$ is the first $L_s$ components and $\hat{q}^1$ is the last $L-L_s$ components.
A similar split can be done for $\hat{H}$, then the mass matrix and the time-stiffness matrix
can be written as
\begin{equation}
  M=\begin{bmatrix}
    M^{00} & M^{01} \\
    M^{10} & M^{11} \\
  \end{bmatrix},
  K_\tau=\begin{bmatrix}
    K_\tau^{00} & K_\tau^{01} \\
    K_\tau^{10} & K_\tau^{11} \\
  \end{bmatrix},
\end{equation}
where the dimensions of sub-matrices $M^{00},M^{01},M^{10},M^{11}$ are
$L_s\times L_s,L_s\times (L-L_s),(L-L_s)\times L_s,(L-L_s)\times (L-L_s)$ respectively,
and it is similar for the sub-matrices of $K_\tau$.
In  \eqref{eq:CG_predictor}, only the last $L-L_s$ components are useful,
and they can be written as
\begin{equation}
  \hat{q}^1=-\hat{M}\hat{h}^1-\hat{M}^0\hat{h}^0,
  \label{eq:CG_iter}
\end{equation}
where $\hat{M}=(K_\tau^{11})^{-1}M^{11},\hat{M}^0=(K_\tau^{11})^{-1}M^{10}$.
We can obtain $\hat{q}^1$ from $\hat{h}^1$ through one iteration using the above equation.
In the continuous Galerkin predictor step, $M$ times iterations of Eq. \eqref{eq:CG_iter}
are adequate for a $M$-th order scheme \cite{Balsara2013}, thus the cost of the iterative part in our scheme is not high.
Once the basis functions are determined, the matrices in  \eqref{eq:CG_iter} are known,
and the whole iterative scheme can be explicitly written down.
We are going to describe the implementation details in the next subsection.

\subsection{Implementation details of  2D third order continuous Galerkin predictor}
This subsection will give the implementation details at third order.
Other cases can be completed similarly and we will provide some difference of the implementation in other cases
in Remark \ref{rem:points}.
Assume that the basis functions in the reference element $[-\frac12,\frac12]$ are orthogonal Legendre polynomials
\begin{equation}
  P_0(\xi)=1,\quad P_1(\xi)=\xi,\quad P_2(\xi)=\xi^2-\dfrac{1}{12},\quad
  P_3(\xi)=\xi^3-\dfrac{3}{20}\xi.
\end{equation}
Thus the solution at time $t^n$ or $\tau=0$ can be represented as a combination of $L_s=6$ basis functions
\begin{align}
  \varphi_h(x,t^n)=&\hat{w}_0P_0(\xi)P_0(\eta)+ \hat{w}_1P_1(\xi)P_0(\eta)+
  \hat{w}_2P_0(\xi)P_1(\eta) \nonumber \\
  &+\hat{w}_3P_2(\xi)P_0(\eta)+ \hat{w}_4P_0(\xi)P_2(\eta)+
  \hat{w}_5P_1(\xi)P_1(\eta).
  \label{eq:spatial_basis}
\end{align}
Take the basis functions in the temporal reference element $[0,1]$ as
\begin{equation}
  Q_0(\tau)=1,\quad Q_1(\tau)=\tau,\quad Q_2(\tau)=\tau^2,\quad
  Q_3(\tau)=\tau^3,
\end{equation}
the first three of which are needed for the third order scheme
while the last basis function is only needed for fourth order schemes.
In order to obtain full third order accuracy in space and time we use a total of $L=10$ basis functions,
and the continuous Galerkin predictor solution can be expressed as
\begin{align}
  q_h(\xi,\eta,\tau)=&\hat{w}_0P_0(\xi)P_0(\eta)Q_0(\tau)+ \hat{w}_1P_1(\xi)P_0(\eta)Q_0(\tau)+
  \hat{w}_2P_0(\xi)P_1(\eta)Q_0(\tau) \nonumber \\
  &+\hat{w}_3P_2(\xi)P_0(\eta)Q_0(\tau)+ \hat{w}_4P_0(\xi)P_2(\eta)Q_0(\tau)+
  \hat{w}_5P_1(\xi)P_1(\eta)Q_0(\tau) \nonumber \\
  &+\hat{q}_6P_0(\xi)P_0(\eta)Q_1(\tau)+ \hat{q}_7P_1(\xi)P_0(\eta)Q_1(\tau)+
  \hat{q}_8P_0(\xi)P_1(\eta)Q_1(\tau) \nonumber \\
  &+\hat{q}_9P_0(\xi)P_0(\eta)Q_2(\tau),
\end{align}
noticing that the first $L_s$ coefficients have been substituted by $\hat{w}_l,l=0,\cdots,L_s-1$ due to the simplification.

Now we can explicitly write down the iterative equation \eqref{eq:CG_iter}
\begin{equation}
    \hat{q}_6 = -\hat{h}_0+\frac{3}{10}\hat{h}_7, ~\hat{q}_7 = -\frac12\hat{h}_6-\frac{3}{5}\hat{h}_7,
    ~\hat{q}_8 = -\hat{h}_1-\frac{2}{3}\hat{h}_8, ~\hat{q}_9 =
    -\hat{h}_2-\frac{2}{3}\hat{h}_9,
\end{equation}
which give one iteration in the predictor step.

At last we have to obtain $\hat{h}$ from $\hat{q}$. The most accurate way is to use the $L^2$ projection,
but it is expensive to do lots of numerical integration. Following \cite{Balsara2013},
we use the nodal approach to determine $\hat{h}$.
Choose $L_n$ points in the reference element, where $L_n$ is equal to or slightly larger than $L$.
In the two-dimensional case at third order, we choose $L_n=13$ nodal points
\begin{align}
  \Big\{
    &(\frac12,0,0),(-\frac12,0,0),(0,\frac12,0),(0,-\frac12,0),
    (\frac12,\frac12,0),(-\frac12,\frac12,0),(\frac12,-\frac12,0),(-\frac12,-\frac12,0),\nonumber \\
    &(\frac12,0,\frac12),(-\frac12,0,\frac12),(0,\frac12,\frac12),(0,-\frac12,\frac12), (0,0,1)
  \Big\}.
\end{align}
Then we can define two $L_n$ component vectors $\bar{u}$ and $\bar{v}$,
which are the nodal values of $\pd{q_h}{\xi}$ and $\pd{q_h}{\eta}$,
and their ordering follows that of nodal points.
Without ambiguity in the text, here we use bar to denote the values at nodal points,
which are not the cell average values defined in the numerical fluxes.
$\bar{u}$ and $\bar{v}$ can be written down explicitly
\begin{align}
  &\bar{u}_0=\hat{q}_1+\hat{q}_3,~\bar{u}_1=\hat{q}_1-\hat{q}_3,
  ~\bar{u}_2=\hat{q}_1+\frac12\hat{q}_5,~\bar{u}_3=\hat{q}_1-\frac12\hat{q}_5,\nonumber \\
  &\bar{u}_4 = \bar{u}_0+\frac12\hat{q}_5,~\bar{u}_5 = \bar{u}_1+\frac12\hat{q}_5,
  ~\bar{u}_6 = \bar{u}_0-\frac12\hat{q}_5,~\bar{u}_7 = \bar{u}_1-\frac12\hat{q}_5,\nonumber \\
  &\bar{u}_8 = \bar{u}_0+\frac12\hat{q}_8,~\bar{u}_9 = \bar{u}_1+\frac12\hat{q}_8,
  ~\bar{u}_{10} = \bar{u}_2+\frac12\hat{q}_8,~\bar{u}_{11} =
  \bar{u}_3+\frac12\hat{q}_8,~\bar{u}_{12} = \hat{q}_1+\hat{q}_8,\nonumber \\
  &\bar{v}_0=\hat{q}_2+\frac12\hat{q}_5,~\bar{v}_1=\hat{q}_2-\frac12\hat{q}_5,
  ~\bar{v}_2=\hat{q}_2+\hat{q}_4,~\bar{v}_3=\hat{q}_2-\hat{q}_4,\nonumber \\
  &\bar{v}_4 = \bar{v}_0+\hat{q}_4,~\bar{v}_5 = \bar{v}_1+\hat{q}_4,
  ~\bar{v}_6 = \bar{v}_0-\hat{q}_4,~\bar{v}_7 = \bar{v}_1-\hat{q}_4,\nonumber \\
  &\bar{v}_8 = \bar{v}_0+\frac12\hat{q}_9,~\bar{v}_9 = \bar{v}_1+\frac12\hat{q}_9,
  ~\bar{v}_{10} = \bar{v}_2+\frac12\hat{q}_9,~\bar{v}_{11} =
  \bar{v}_3+\frac12\hat{q}_9,~\bar{v}_{12} = \hat{q}_2+\hat{q}_9,
\end{align}
also noticing that the first eight elements of $\bar{u}$ and $\bar{v}$ only depend on the initial condition,
thus they only have to be calculated once in the predictor step.
Once we obtain $\bar{u}$ and $\bar{v}$, we can get the values of $h$ at each nodal points respectively,
denoted by $\bar{h}$.
The above formulae give a transcription from modal space to nodal space.
Then we may express the transcription from nodal space to modal space explicitly,
where $r_1,r_2$ are temporary variables to save cost in the calculation
\begin{align}
  &\hat{h}_1=\bar{h}_0-\bar{h}_1,~
  \hat{h}_2=\bar{h}_2-\bar{h}_3,~
  \hat{h}_5=2(\bar{h}_4-\bar{h}_5-\hat{h}_1),~
  \hat{h}_3=4(\bar{h}_4-\bar{h}_2)-2\hat{h}_1-\hat{h}_5,\nonumber\\
  &\hat{h}_4=4(\bar{h}_4-\bar{h}_0)-2\hat{h}_2-\hat{h}_5,~
  \hat{h}_0=\frac18(\bar{h}_0+\bar{h}_1+\bar{h}_2+\bar{h}_3+\bar{h}_4+\bar{h}_5+\bar{h}_6+\bar{h}_7-\frac56(\hat{h}_3+\hat{h}_4)),\nonumber\\
  &\hat{h}_8=2(\bar{h}_8-\bar{h}_9-\bar{h}_0+\bar{h}_1),~
  \hat{h}_9=2(\bar{h}_{10}-\bar{h}_{11}-\bar{h}_2+\bar{h}_3),\nonumber\\
  &r_1=\bar{h}_9+\bar{h}_9-\bar{h}_0-\bar{h}_1,~
  r_2=\bar{h}_{12}-\hat{h}_0+\frac{1}{12}(\hat{h}_3+\hat{h}_4),~
  \hat{h}_7=2(r_2-r_1),~
  \hat{h}_6=r_2-\hat{h}_7.
\end{align}
Similarly, we only have to compute the first six coefficients once,
while the last four need to be determined in the iterative predictor step.
So far we have provided all the implementation details of the continuous Galerkin predictor step
in the two-dimensional case at third order.

\begin{remark}\label{rem:points}
  For the strategy to choose the points in the reference element, we refer the readers to \cite{Balsara2013}.
  In this paper, the points we choose in the one-dimensional case are
  \begin{enumerate}
    \item second order: $\left\{(\frac12,0),(-\frac12,0),(0,1) \right\}$,
    \item third order:
      $\left\{(0,0),(\frac12,0),(-\frac12,0),(\frac12,\frac12),(-\frac12,\frac12),(0,1)
      \right\}$,
    \item fourth order:
      $\left\{(0,0),(\frac12,0),(-\frac12,0),(\frac14,0),(-\frac14,0),
      (0,\frac13),(\frac12,\frac13),(-\frac12,\frac13),
    (\frac12,\frac23),(-\frac12,\frac23),(0,1) \right\}$.
  \end{enumerate}
  The points we choose in the two-dimensional case are
  \begin{enumerate}
    \item second order:
  $\left\{ (\frac12,0,0),(-\frac12,0,0),(0,\frac12,0),(0,-\frac12,0),(0,0,1) \right\}$,
    \item fourth order:
      \begin{align*}
        \big\{
          &(0,0,0),(\frac12,0,0),(-\frac12,0,0),(0,\frac12,0),(0,-\frac12,0),
          (\frac12,\frac12,0),(-\frac12,\frac12,0),(\frac12,-\frac12,0),(-\frac12,-\frac12,0),\\
          &(\frac14,0,0),(-\frac14,0,0),(0,\frac14,0),(0,-\frac14,0),
          (0,0,\frac13),(\frac12,\frac12,\frac13),(-\frac12,\frac12,\frac13),(\frac12,-\frac12,\frac13),(-\frac12,-\frac12,\frac13),\\
          &(\frac12,0,\frac23),(-\frac12,0,\frac23),(0,\frac12,\frac23),(0,-\frac12,\frac23),
          (0,0,1)
        \big\}.
      \end{align*}
  \end{enumerate}
  Once those points are determined, the transcription between modal space and nodal space can be obtained by the similar procedure described for two-dimensional case at third order.
\end{remark}

\section{Calculation of the volume integral and the numerical fluxes in the ADER-DG scheme}\label{sec:flux}
In the RKDG scheme \cite{Cheng2014}, we have to calculate the volume integral of Hamiltonian,
and the numerical fluxes (in two-dimensional case or above) by numerical integration,
so that we have to compute Roe speed and the parameter at every integration point at cell interfaces,
which is expensive.
If we use the same way to accomplish them in the ADER-DG scheme,
the computational complexity will be much larger, because one more dimension will appear in the ADER-DG scheme.
For the volume integral, thanks to the spacetime representation of the Hamiltonian in the spacetime control volume,
it can be expressed explicitly.
For the numerical fluxes, we use a substantially simpler strategy presented in \cite{Dumbser2007} to calculate
the numerical fluxes in our ADER-DG scheme.

\subsection{The explicit formulae of the volume integral}
Because we have saved all $\hat{h}$ in each cell, i.e. the coefficients of the basis functions,
and the test functions $\phi_m$ in \eqref{eq:ADER-DG2D} are transformed from the six basis functions
in \eqref{eq:spatial_basis} to the computational cell $I_{i,j}$,
we can write down the spacetime integral of the Hamiltonian in the ADER-DG scheme \eqref{eq:ADER-DG2D} immediately
\begin{align}
  &\intt\int_{I_{i,j}} \phi_0H(\partial_x q_h,\partial_y q_h,x,y)\dd x\dd y\dd t
  =\Delta x_i\Delta y_j(\hat{h}_0+\frac12\hat{h}_6+\frac13\hat{h}_7),\nonumber\\
  &\intt\int_{I_{i,j}} \phi_1H(\partial_x q_h,\partial_y q_h,x,y)\dd x\dd y\dd t
  =\Delta x_i\Delta y_j(\frac{1}{12}\hat{h}_1+\frac{1}{24}\hat{h}_8),\nonumber\\
  &\intt\int_{I_{i,j}} \phi_2H(\partial_x q_h,\partial_y q_h,x,y)\dd x\dd y\dd t
  =\Delta x_i\Delta y_j(\frac{1}{12}\hat{h}_2+\frac{1}{24}\hat{h}_9),\nonumber\\
  &\intt\int_{I_{i,j}} \phi_3H(\partial_x q_h,\partial_y q_h,x,y)\dd x\dd y\dd t
  =\frac{\Delta x_i\Delta y_j}{180}\hat{h}_3,\nonumber\\
  &\intt\int_{I_{i,j}} \phi_4H(\partial_x q_h,\partial_y q_h,x,y)\dd x\dd y\dd t
  =\frac{\Delta x_i\Delta y_j}{180}\hat{h}_4,\nonumber\\
  &\intt\int_{I_{i,j}} \phi_5H(\partial_x q_h,\partial_y q_h,x,y)\dd x\dd y\dd t
  =\frac{\Delta x_i\Delta y_j}{144}\hat{h}_5.
\end{align}

\subsection{The explicit formulae of the numerical fluxes}
In \cite{Dumbser2007}, the central idea consists of freezing the wave speeds
to equal their values evaluated at the spacetime barycenters of the face under consideration.
For the numerical fluxes  considered here, we would freeze the Roe speed $\tilde H_{q_h}$
and the parameter $S_{q_h}$ to   their values at the spacetime barycenters,
then the remaining parts can be explicitly calculated.
Now we take the face
$x_{i+\frac12}\times[y_{j-\frac12},y_{j+\frac12}]\times[t^n,t^{n+1}]$ for example,
whose left neighbor cell   is $I_{i,j}$, and the  right is $I_{i+1,j}$.
In order to calculate the Roe speed $\tilde H_{1,q_h}$ and the parameter $S_{1,q_h}$
at the spacetime barycenter $(x_{i+\frac12},y_{j},t^{n+\frac12})$,
we need to give the left and right limit values of the partial derivatives
at the point
\begin{align}
  &u_L=\hat{q}_{L,1}+\hat{q}_{L,3}+\frac12\hat{q}_{L,8},
  ~u_R=\hat{q}_{R,1}-\hat{q}_{R,3}+\frac12\hat{q}_{R,8},\nonumber\\
  &v_L=\hat{q}_{L,2}+\frac12\hat{q}_{L,5}+\frac12\hat{q}_{L,9},
  ~v_R=\hat{q}_{R,2}-\frac12\hat{q}_{R,5}+\frac12\hat{q}_{R,9},
  \label{eq:para}
\end{align}
where the subscripts $L$ and $R$ denote values in the left and right side respectively,
and $u,v$ denote the partial derivatives $(q_h)_x$ and $(q_h)_y$ respectively.
Then we can obtain $\tilde H_{1,q_h}$ and $S_{1,q_h}$ by using the definition \eqref{eq:para}
at $(x_{i+\frac12},y_{j},t^{n+\frac12})$.
Now freeze $\tilde H_{1,q_h}$ and $S_{1,q_h}$ in the numerical fluxes terms in the ADER-DG scheme \eqref{eq:ADER-DG2D},
and introduce three temporary variables as
\begin{align*}
&\lambda_1=\min(\tilde{H}_{1,{q_h}}(x_\ir,y_j,t^{n+\frac12}),0),\\
&\lambda_2=\max(\tilde{H}_{1,{q_h}}(x_\ir,y_j,t^{n+\frac12}),0),\\
&\lambda_3=S_{1,q_h}(x_\ir,y_j,t^{n+\frac12})-\abs{\tilde{H_{1,q_h}}(x_\ir,y_j,t^{n+\frac12})}.
\end{align*}
Denote the basis functions in the left cell $I_{i,j}$ and in the right cell $I_{i+1,j}$
by $\phi_{L,m}$ and $\phi_{R,m}$ respectively.
For the first basis function $\phi_{L,0}=\phi_{R,0}=1$, we have
\begin{align}
  &\varphi_L=\hat{q}_{L,0}+\frac12(\hat{q}_{L,1}+\hat{q}_{L,6})+\frac16\hat{q}_{L,3}
  +\frac13\hat{q}_{L,7}+\frac14\hat{q}_{L,8},\nonumber \\
  &\varphi_R=\hat{q}_{R,0}+\frac12(-\hat{q}_{R,1}+\hat{q}_{R,6})+\frac16\hat{q}_{R,3}
  +\frac13\hat{q}_{R,7}-\frac14\hat{q}_{R,8},\nonumber \\
  &u_L=\hat{q}_{L,1}+\hat{q}_{L,3}+\frac12\hat{q}_{L,8},\nonumber \\
  &u_R=\hat{q}_{R,1}-\hat{q}_{R,3}+\frac12\hat{q}_{R,8},
\end{align}
and the numerical fluxes terms in \eqref{eq:ADER-DG2D} are
\begin{align}
  &\intt\int_{K_j}\min(\tilde{H}_{1,{q_h}}(x_\ir,y,t),0)[{q_h}](x_\ir,y,t)\phi_{L,0}(x_\ir^-,y)\dd y\dd t
  \approx \lambda_1(\varphi_R-\varphi_L)\Delta t\Delta y_j, \nonumber \\
  &\intt\int_{K_j}\max(\tilde{H}_{1,{q_h}}(x_\ir,y,t),0)[{q_h}](x_\ir,y,t)\phi_{R,0}(x_\il^+,y)\dd y\dd t
  \approx \lambda_2(\varphi_R-\varphi_L)\Delta t\Delta y_j, \nonumber \\
  &C\Delta x_i\intt\int_{K_j}\left(
  S_{1,{q_h}}(x_\ir,y,t)-\abs{\tilde{H}_{1,{q_h}}(x_\ir,y,t)}\right)
  [(q_h)_x](x_\ir,y)\phi_{L,0}(x_\ir^-,y)\dd y\dd t \nonumber\\
  &\approx C\lambda_3(u_R-u_L)\Delta t\Delta y_j, \nonumber \\
  &C\Delta x_i\intt\int_{K_j}\left(
  S_{1,{q_h}}(x_\ir,y,t)-\abs{\tilde{H}_{1,{q_h}}(x_\ir,y,t)}\right)
  [(q_h)_x](x_\ir,y)\phi_{R,0}(x_\il^+,y)\dd y\dd t \nonumber\\
  &\approx C\lambda_3(u_R-u_L)\Delta t\Delta y_j,
  \label{eq:bd_integral}
\end{align}
where the first and the third integrations are the contributions to the left cell,
and the second and the fourth integrations are the contributions to the right cell.
For the second basis function transformed from $\phi_{1}=P_1(\xi(x))P_0(\eta(y))$ and the fourth basis function
transformed from $\phi_{3}=P_2(\xi(x))P_0(\eta(y))$, the corresponding integrations are
$$ \frac{1}{2}\lambda_1(\varphi_R-\varphi_L)\Delta t\Delta y_j,
  ~-\frac{1}{2}\lambda_2(\varphi_R-\varphi_L)\Delta t\Delta y_j,
  ~\frac{1}{2}C\lambda_3(u_R-u_L)\Delta t\Delta y_j,
  ~-\frac{1}{2}C\lambda_3(u_R-u_L)\Delta t\Delta y_j,$$
and
$$ \frac{1}{6}\lambda_1(\varphi_R-\varphi_L)\Delta t\Delta y_j,
  ~\frac{1}{6}\lambda_2(\varphi_R-\varphi_L)\Delta t\Delta y_j,
  ~\frac{1}{6}C\lambda_3(u_R-u_L)\Delta t\Delta y_j,
  ~\frac{1}{6}C\lambda_3(u_R-u_L)\Delta t\Delta y_j.$$
We can clearly see from the above formulae that, for the second and the fourth basis functions,
the numerical fluxes terms are just a scaling of the corresponding numerical fluxes terms for the first basis function,
thus the computational costs can be reduced greatly.
Similarly, if denoting
\begin{align}
  &\varphi_L= \frac{1}{24}(2\hat{q}_{L,2}+\hat{q}_{L,5}+\hat{q}_{L,9}),
  ~\varphi_R= \frac{1}{24}(2\hat{q}_{R,2}-\hat{q}_{R,5}+\hat{q}_{R,9}),\nonumber \\
  &u_L=\frac{1}{12}\hat{q}_{L,5},
  ~u_R=\frac{1}{12}\hat{q}_{R,5},
\end{align}
then for the third basis function transformed from $\phi_{2}=P_0(\xi(x))P_1(\eta(y))$
and the sixth basis function transformed from $\phi_{5}=P_1(\xi(x))P_1(\eta(y))$,
the four numerical fluxes terms in \eqref{eq:ADER-DG2D} corresponding to \eqref{eq:bd_integral} are
$$\lambda_1(\varphi_R-\varphi_L)\Delta t\Delta y_j,
  ~\lambda_2(\varphi_R-\varphi_L)\Delta t\Delta y_j,
  ~C\lambda_3(u_R-u_L)\Delta t\Delta y_j,
  ~C\lambda_3(u_R-u_L)\Delta t\Delta y_j,$$
and
$$\frac{1}{2}\lambda_1(\varphi_R-\varphi_L)\Delta t\Delta y_j,
  ~-\frac{1}{2}\lambda_2(\varphi_R-\varphi_L)\Delta t\Delta y_j,
  ~\frac{1}{2}C\lambda_3(u_R-u_L)\Delta t\Delta y_j,
  ~-\frac{1}{2}C\lambda_3(u_R-u_L)\Delta t\Delta y_j,$$
respectively.
If denoting
\begin{align}
  &\varphi_L= \frac{1}{180}\hat{q}_{L,4},
  ~\varphi_R= \frac{1}{180}\hat{q}_{R,4},
  ~u_L=0,
  ~u_R=0,
\end{align}
then for the fifth basis function transformed from $\phi_{4}=P_0(\xi(x))P_2(\eta(y))$,
the four corresponding numerical fluxes terms are
$$\lambda_1(\varphi_R-\varphi_L)\Delta t\Delta y_j,
  ~\lambda_2(\varphi_R-\varphi_L)\Delta t\Delta y_j,
  ~C\lambda_3(u_R-u_L)\Delta t\Delta y_j,
  ~C\lambda_3(u_R-u_L)\Delta t\Delta y_j.$$

We have explicitly given the volume integral and the numerical fluxes terms in
the ADER-DG scheme \eqref{eq:ADER-DG2D} in the two-dimensional case at third order.
In the next subsection, we would like to compare the computational complexity
of the ADER-DG scheme to that of the RKDG scheme.

\subsection{Comparison of the computational complexity between the ADER-DG scheme and the RKDG scheme}\label{subsec:CPU}
This section give a comparison of the computational complexities
of the ADER-DG scheme and the RKDG scheme.
As an example, we consider them in the two-dimensional case at third order,
and are going to count the number of operations needed
in the evolution procedure at a time step for one cell.
Four types of basic operations, i.e., addition, subtraction, multiplication and division,
are all treated as one operation.
And we regard one calculation of the Hamiltonian as one operation.
The main part is evolving  the DG solutions in each cell,
so the part for calculating the time step is negligible.

In the ADER-DG scheme, we need 79 operations to calculate the partial derivatives
$u$ and $v$ and 31 operations to obtain the value of Hamiltonian,
120 operations to accomplish the transcription from nodal space to modal space,
57 operations to perform three iterations, 13 operations
to calculate the volume integral of Hamiltonian. For each face, we need 28 operations
to compute the Roe speed and the parameters, 84 operations for the numerical fluxes terms.
Because there are four faces for a cell and each face is shared by two cells,
the operations on the faces of one cell should be doubled.
At one time step, we need 36 operations to evolve the DG solution in a cell,
 so 560 operations are needed to update one cell in a time step.

In the RKDG scheme, we need to use numerical integration.
We use three points Gauss-Legendre integration on an edge,
and nine points Gauss-Legendre integration in the tensor product form for the volume integral.
Because the values of each basis functions at each integration points in the reference element
will be used many times, we compute them once and save them. We need 11 and 12 operations to compute
the value and the partial derivatives of $\varphi$ respectively.
Thus we need 486 operations to calculate the volume integral of Hamiltonian once in the RKDG scheme, and on each edge,
 222 operations to compute the Roe speed and the parameter $S$, 189 operations to compute the numerical fluxes.
 At each sub-step, we need 36 operations to evolve the DG solution,
 thus 1344 operations are needed in a sub-step and 4032 operations in total.

At second order, the number of the operations for the ADER-DG scheme and the RKDG scheme are
180 and 868 respectively, and at fourth order, they are 2358 and 12944 respectively.

The above analysis shows that  the ADER-DG scheme has much less computational complexity
when the solution is evolved in a cell at a time step,
which is about $20.7\%,13.9\%,18.2\%$ of the RKDG scheme at second, third and fourth order, respectively.
The reason is that the ADER-DG scheme is a one-step one-stage  scheme
and we use a cheap way to calculate the volume integral and
the numerical fluxes terms in the ADER-DG schemes.
Further, we will show that our schemes can achieve the designed order of accuracy in the numerical experiments,
and the comparison of the CPU times will be also recorded to validate the efficiency of the ADER-DG scheme.

\section{Numerical results}\label{sec:numer}
This section will provide some numerical experiments in one- and two-dimensions.
In the two-dimensional experiments, we use $N\times N$ uniform meshes with $\Delta x=\Delta y$.
The time stepsize  is chosen as $\Delta t = {\text{CFL}\Delta x}/{\alpha}$,
where $\alpha=\max{\abs{H_1}}$ for one-dimensional cases
and $\alpha=\max\{\abs{H_1},\abs{H_2}\}$ for two-dimensional cases.

\subsection{One-dimensional results}
\begin{example}
{\rm \label{ex:variable_coeff_sinx}
	We solve the following linear problem \cite{Cheng2007} with a smooth variable
  coefficient
	\begin{equation*}
	\varphi_t+\sin(x)\varphi_x=0,\quad 0\leqslant x\leqslant 2\pi.
	\end{equation*}
	The initial condition is $\varphi(x,0)=\sin(x)$, and the periodic boundary condition is specified.
	The exact solution is $\varphi(x,t)=\sin(2\arctan(e^{-t}\tan(\frac x2)))$.
}
\end{example}

The numerical errors and the orders of convergence at $t=1$ are presented in Table \ref{tab:sinx1D}.
We can see that the ADER-DG scheme can achieve $(k+1)$-th order accuracy for $P^k$ polynomials.

\begin{table}[!ht]
  \centering
  \begin{tabular}{c|cccccc}
    \hline $N$ & $\ell^2$ error & Order & $\ell^1$ error & Order & $\ell^\infty$ error & Order \\
    \hline \multicolumn{7}{c}{$k=1,\text{CFL}=0.15$} \\ \hline
 20 & 1.322e-02 &  -   & 1.884e-02 &  -   & 1.966e-02 &  -   \\
 40 & 3.620e-03 & 1.87 & 5.096e-03 & 1.89 & 5.318e-03 & 1.89 \\
 80 & 9.576e-04 & 1.92 & 1.392e-03 & 1.87 & 1.524e-03 & 1.80 \\
160 & 2.437e-04 & 1.97 & 3.458e-04 & 2.01 & 3.891e-04 & 1.97 \\
320 & 6.157e-05 & 1.98 & 8.612e-05 & 2.01 & 9.859e-05 & 1.98 \\
640 & 1.536e-05 & 2.00 & 2.162e-05 & 1.99 & 2.500e-05 & 1.98 \\
    \hline \multicolumn{7}{c}{$k=2,\text{CFL}=0.10$} \\ \hline
 20 & 1.060e-03 &  -   & 1.460e-03 &  -   & 1.761e-03 &  -   \\
 40 & 1.391e-04 & 2.93 & 2.022e-04 & 2.85 & 1.976e-04 & 3.16 \\
 80 & 2.033e-05 & 2.77 & 2.781e-05 & 2.86 & 3.556e-05 & 2.47 \\
160 & 2.868e-06 & 2.83 & 3.778e-06 & 2.88 & 5.535e-06 & 2.68 \\
320 & 3.927e-07 & 2.87 & 5.063e-07 & 2.90 & 7.534e-07 & 2.88 \\
640 & 5.230e-08 & 2.91 & 6.658e-08 & 2.93 & 9.928e-08 & 2.92 \\
    \hline \multicolumn{7}{c}{$k=3,\text{CFL}=0.05$} \\ \hline
 20 & 7.609e-05 &  -   & 1.278e-04 &  -   & 9.354e-05 &  -   \\
 40 & 8.493e-06 & 3.16 & 1.034e-05 & 3.63 & 2.101e-05 & 2.15 \\
 80 & 6.436e-07 & 3.72 & 7.851e-07 & 3.72 & 1.307e-06 & 4.01 \\
160 & 4.450e-08 & 3.85 & 5.410e-08 & 3.86 & 7.630e-08 & 4.10 \\
320 & 2.939e-09 & 3.92 & 3.557e-09 & 3.93 & 4.779e-09 & 4.00 \\
640 & 1.890e-10 & 3.96 & 2.278e-10 & 3.96 & 3.141e-10 & 3.93 \\
    \hline
  \end{tabular}
  \caption{Errors and orders of convergence for Example \ref{ex:variable_coeff_sinx},~$t=1$}
  \label{tab:sinx1D}
\end{table}

\begin{example}
  {\rm \label{ex:sign1D}
  	We solve the following linear problem  \cite{Cheng2007}
  	\begin{equation*}
  	\varphi_t+\text{sign}(\cos(x))\varphi_x=0,\quad 0\leqslant x\leqslant 2\pi,
  	\end{equation*}
  	with initial condition $\varphi(x,0)=\sin(x)$, and periodic boundary condition.
Obviously,  the variable coefficient is not smooth.
  }
\end{example}

\begin{figure}[!ht]
  \centering
  \subfigure[$P^2$]{
    \includegraphics[width=0.45\textwidth]{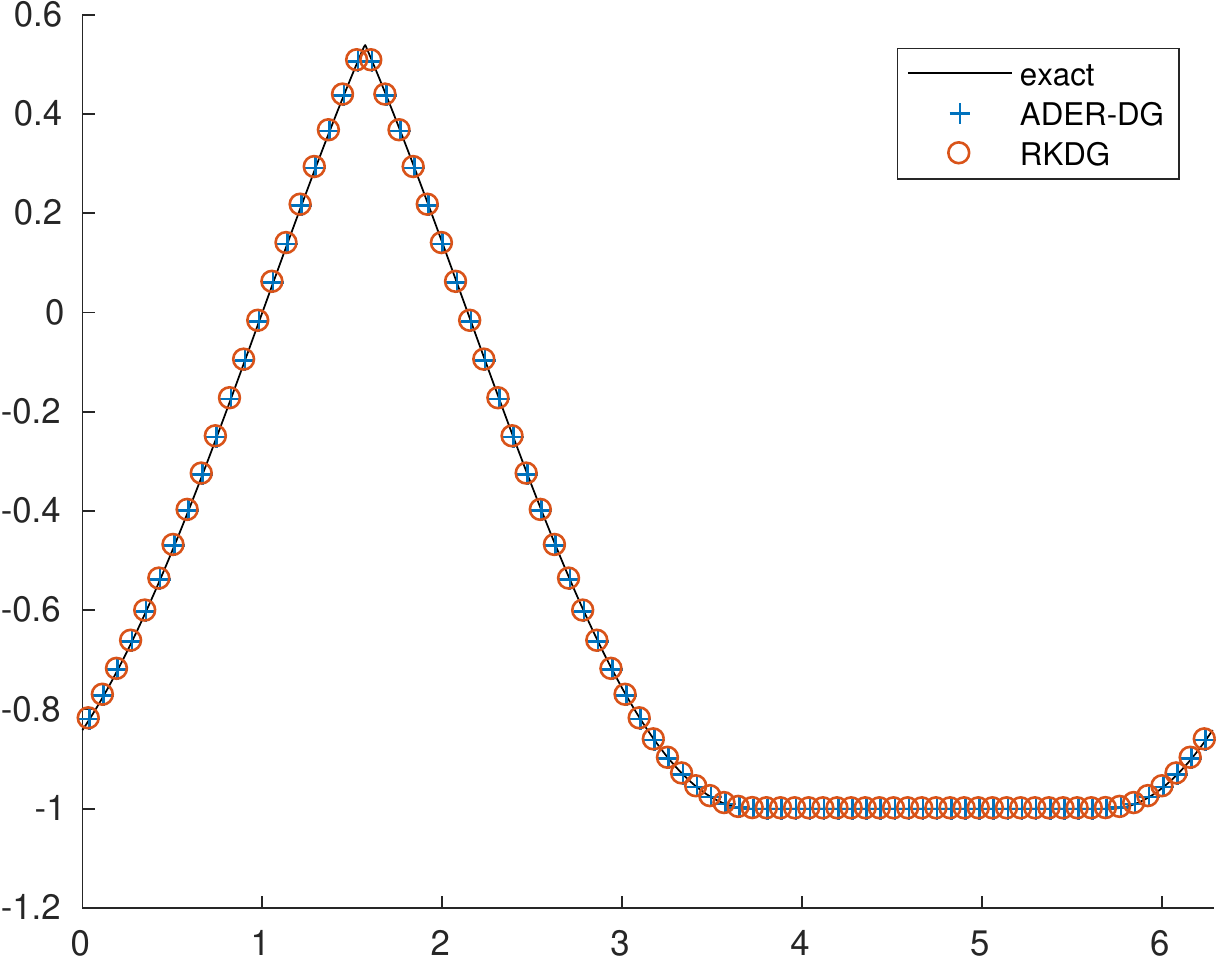}
  }
  \subfigure[$P^3$]{
    \includegraphics[width=0.45\textwidth]{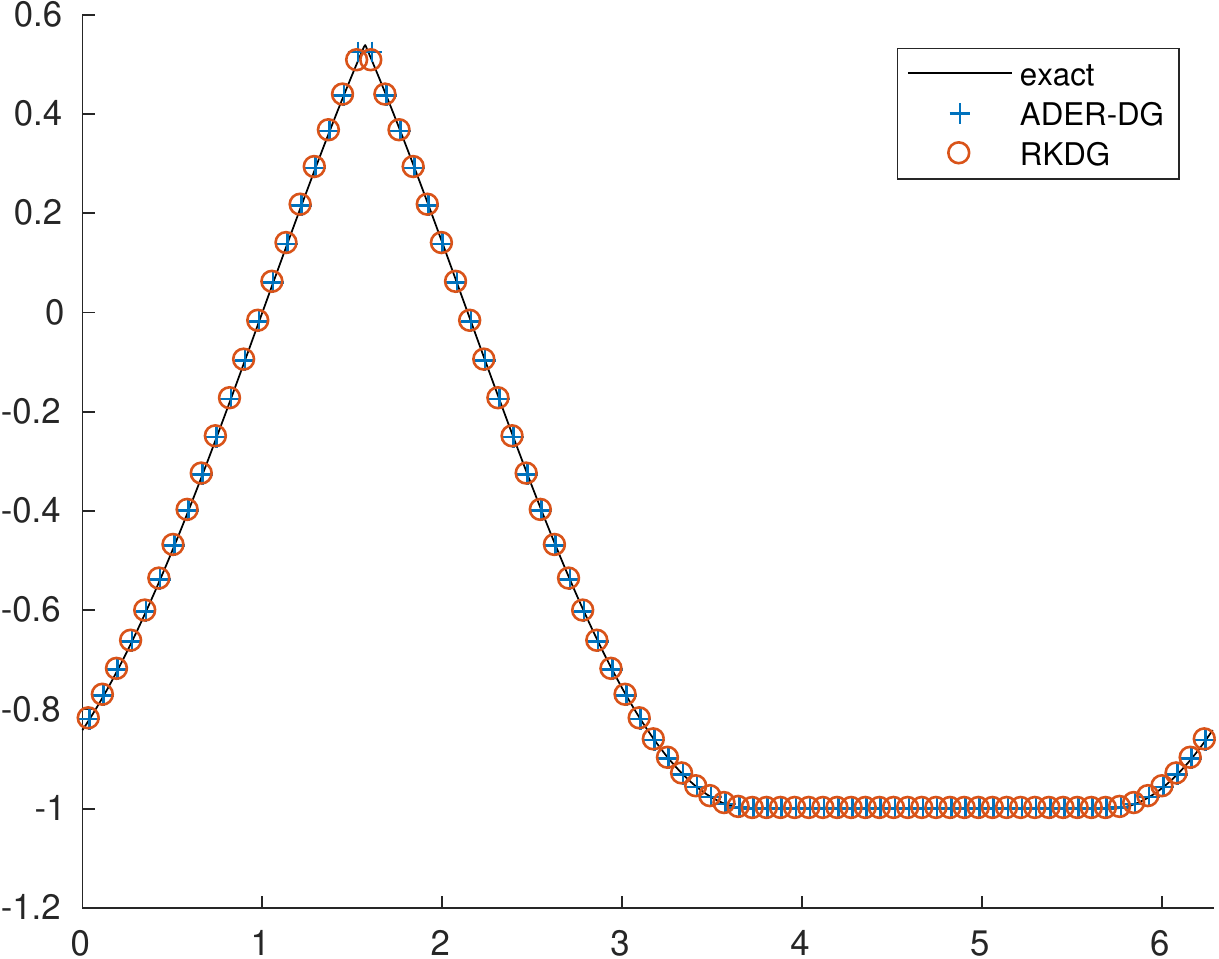}
  }
  \caption{Example \ref{ex:sign1D},~$N=80$. Left: $P^2$, right: $P^3$. }
  \label{fig:sign1D}
\end{figure}

In the viscosity solution, there is a shock forming in $\varphi_x$ at $x=\frac{\pi}{2}$, and a rarefaction wave at $x=\frac{3\pi}{2}$,
thus the numerical errors are only calculated in the smooth region $[0,1]\cup [2,3.4]\cup [6,2\pi]$.
The errors and the orders of convergence at $t=1$ are presented in Table \ref{tab:sign1D}.
From the table, we can observe that our schemes can achieve $(k+1)$-th order accuracy for $P^k$ polynomials in the smooth region.
The results obtained with $P^2$ and $P^3$ ADER-DG scheme and   $N=80$
are also shown in Fig. \ref{fig:sign1D}.
The ADER-DG scheme can converge to the viscosity solution.

\begin{table}[!ht]
  \centering
  \begin{tabular}{c|cccccc}
    \hline $N$ & $\ell^2$ error & Order & $\ell^1$ error & Order & $\ell^\infty$ error & Order \\
    \hline \multicolumn{7}{c}{$k=1,\text{CFL}=0.10$} \\ \hline
 20 & 9.156e-03 &  -   & 1.276e-02 &  -   & 5.403e-03 &  -   \\
 40 & 2.206e-03 & 2.05 & 3.023e-03 & 2.08 & 1.515e-03 & 1.83 \\
 80 & 3.920e-04 & 2.49 & 5.679e-04 & 2.41 & 4.227e-04 & 1.84 \\
160 & 9.269e-05 & 2.08 & 1.355e-04 & 2.07 & 5.170e-05 & 3.03 \\
320 & 2.111e-05 & 2.13 & 3.074e-05 & 2.14 & 1.341e-05 & 1.95 \\
640 & 5.127e-06 & 2.04 & 7.532e-06 & 2.03 & 2.425e-06 & 2.47 \\
    \hline \multicolumn{7}{c}{$k=2,\text{CFL}=0.10$} \\ \hline
 20 & 3.419e-04 &  -   & 4.731e-04 &  -   & 3.417e-04 &  -   \\
 40 & 3.845e-05 & 3.15 & 5.057e-05 & 3.23 & 4.548e-05 & 2.91 \\
 80 & 4.185e-06 & 3.20 & 5.476e-06 & 3.21 & 5.279e-06 & 3.11 \\
160 & 4.920e-07 & 3.09 & 6.129e-07 & 3.16 & 6.623e-07 & 2.99 \\
320 & 6.144e-08 & 3.00 & 7.668e-08 & 3.00 & 8.205e-08 & 3.01 \\
640 & 7.703e-09 & 3.00 & 9.629e-09 & 2.99 & 1.040e-08 & 2.98 \\
    \hline \multicolumn{7}{c}{$k=3,\text{CFL}=0.03$} \\ \hline
 20 & 3.168e-05 &  -   & 2.887e-05 &  -   & 5.832e-05 &  -   \\
 40 & 2.891e-06 & 3.45 & 2.017e-06 & 3.84 & 5.459e-06 & 3.42 \\
 80 & 9.858e-08 & 4.87 & 5.897e-08 & 5.10 & 4.190e-07 & 3.70 \\
160 & 1.298e-09 & 6.25 & 1.268e-09 & 5.54 & 3.669e-09 & 6.84 \\
320 & 5.215e-11 & 4.64 & 6.367e-11 & 4.32 & 3.469e-11 & 6.72 \\
640 & 3.250e-12 & 4.00 & 3.970e-12 & 4.00 & 2.209e-12 & 3.97 \\
    \hline
  \end{tabular}
  \caption{Errors and orders of convergence for Example \ref{ex:sign1D},~$t=1$}
  \label{tab:sign1D}
\end{table}

\begin{example}
 {\rm \label{ex:burgers1D}
	We solve one-dimensional Burgers' equation
	\begin{equation*}
	\varphi_t+\dfrac{(\varphi_x+1)^2}{2}=0,\quad -1\leqslant x\leqslant 1,
	\end{equation*}
	with smooth initial condition, initial condition $\varphi(x,0)=-\cos(\pi x)$, and periodic boundary condition.
}
\end{example}

\begin{figure}[!ht]
  \centering
  \subfigure[$P^2$]{
    \includegraphics[width=0.45\textwidth]{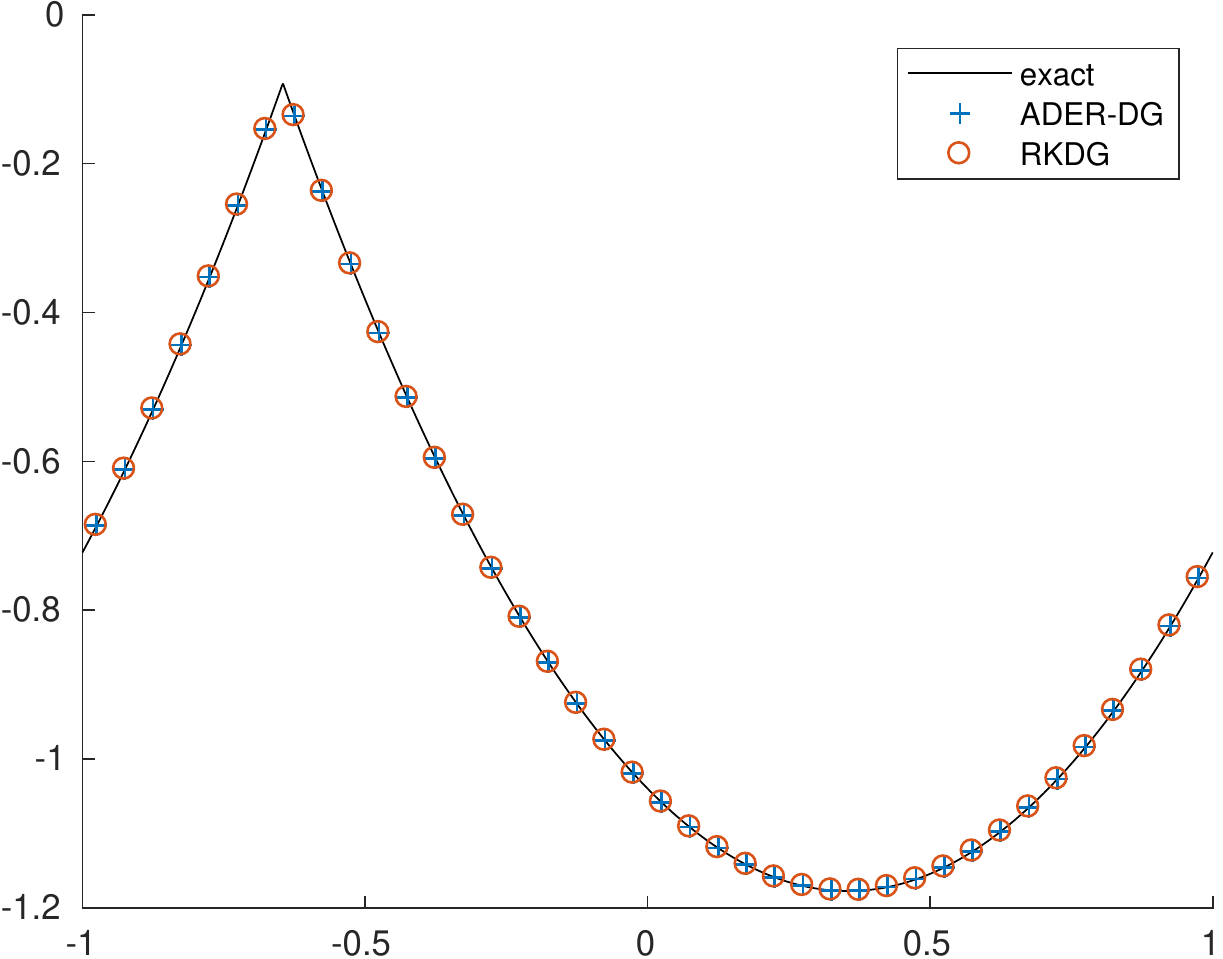}
  }
  \subfigure[$P^3$]{
    \includegraphics[width=0.45\textwidth]{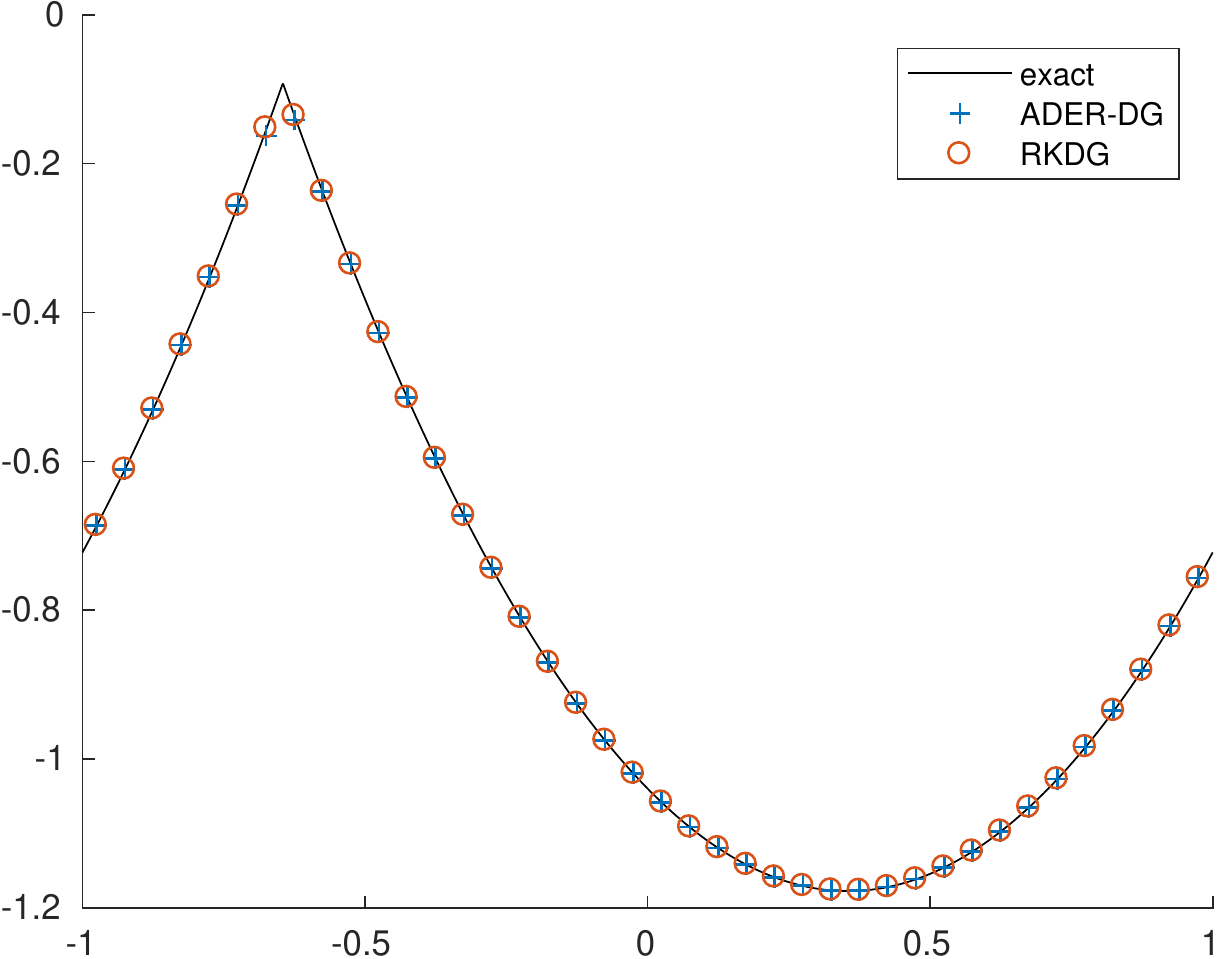}
  }
  \caption{Example \ref{ex:burgers1D},~$N=40$.}
  \label{fig:burgers1D}
\end{figure}

We compute the solution up to $t=\dfrac{0.5}{\pi^2}$.
At this time, the solution is still smooth.
We provide the errors and the orders of convergence in Table \ref{tab:burgers1D}.
Our scheme can achieve the designed order of accuracy in this example.
We also compute the solution up to $t=\dfrac{1.5}{\pi^2}$, there will be a shock in $\varphi_x$.
In Fig. \ref{fig:burgers1D}, we show the results obtained with $P^2$ and $P^3$ ADER-DG scheme with $N=40$.
From the figures, we can see that the ADER-DG scheme can give good results.

\begin{table}[!ht]
  \centering
  \begin{tabular}{c|cccccc}
    \hline $N$ & $\ell^2$ error & Order & $\ell^1$ error & Order & $\ell^\infty$ error & Order \\
    \hline \multicolumn{7}{c}{$k=1,\text{CFL}=0.15$} \\ \hline
 10 & 1.307e-02 &  -   & 1.567e-02 &  -   & 2.185e-02 &  -   \\
 20 & 4.242e-03 & 1.62 & 4.323e-03 & 1.86 & 8.331e-03 & 1.39 \\
 40 & 8.865e-04 & 2.26 & 9.420e-04 & 2.20 & 1.986e-03 & 2.07 \\
 80 & 2.007e-04 & 2.14 & 2.250e-04 & 2.07 & 5.241e-04 & 1.92 \\
160 & 4.850e-05 & 2.05 & 5.566e-05 & 2.02 & 1.534e-04 & 1.77 \\
320 & 1.227e-05 & 1.98 & 1.406e-05 & 1.99 & 4.123e-05 & 1.90 \\
    \hline \multicolumn{7}{c}{$k=2,\text{CFL}=0.10$} \\ \hline
 10 & 1.139e-03 &  -   & 1.244e-03 &  -   & 2.522e-03 &  -   \\
 20 & 1.426e-04 & 3.00 & 1.452e-04 & 3.10 & 3.850e-04 & 2.71 \\
 40 & 2.034e-05 & 2.81 & 2.013e-05 & 2.85 & 5.262e-05 & 2.87 \\
 80 & 2.796e-06 & 2.86 & 2.675e-06 & 2.91 & 7.144e-06 & 2.88 \\
160 & 3.719e-07 & 2.91 & 3.469e-07 & 2.95 & 1.298e-06 & 2.46 \\
320 & 4.843e-08 & 2.94 & 4.456e-08 & 2.96 & 1.924e-07 & 2.75 \\
    \hline \multicolumn{7}{c}{$k=3,\text{CFL}=0.05$} \\ \hline
 10 & 1.320e-04 &  -   & 1.272e-04 &  -   & 2.718e-04 &  -   \\
 20 & 9.644e-06 & 3.78 & 8.445e-06 & 3.91 & 3.740e-05 & 2.86 \\
 40 & 7.291e-07 & 3.73 & 5.760e-07 & 3.87 & 3.211e-06 & 3.54 \\
 80 & 4.937e-08 & 3.88 & 3.793e-08 & 3.92 & 2.347e-07 & 3.77 \\
160 & 3.231e-09 & 3.93 & 2.438e-09 & 3.96 & 1.549e-08 & 3.92 \\
320 & 2.078e-10 & 3.96 & 1.544e-10 & 3.98 & 9.509e-10 & 4.03 \\
    \hline
  \end{tabular}
  \caption{Errors and orders of convergence for Example
  \ref{ex:burgers1D},~$t=\dfrac{0.5}{\pi^2}$}
  \label{tab:burgers1D}
\end{table}

\begin{example}
{\rm \label{ex:burgers_nonsmooth}
	We solve one-dimensional Burgers' equation  \cite{Cheng2014},
	\begin{equation*}
	\varphi_t+\dfrac{\varphi_x^2}{2}=0,\quad 0\leqslant x\leqslant 2\pi,
	\end{equation*}
	with unsmooth initial condition, initial condition $\varphi(x,0)=\abs{x-\pi}$, and periodic boundary condition.
}
\end{example}

\begin{figure}[!ht]
  \centering
  \subfigure[$P^2,\text{CFL}=0.10$]{
    \includegraphics[width=0.45\textwidth]{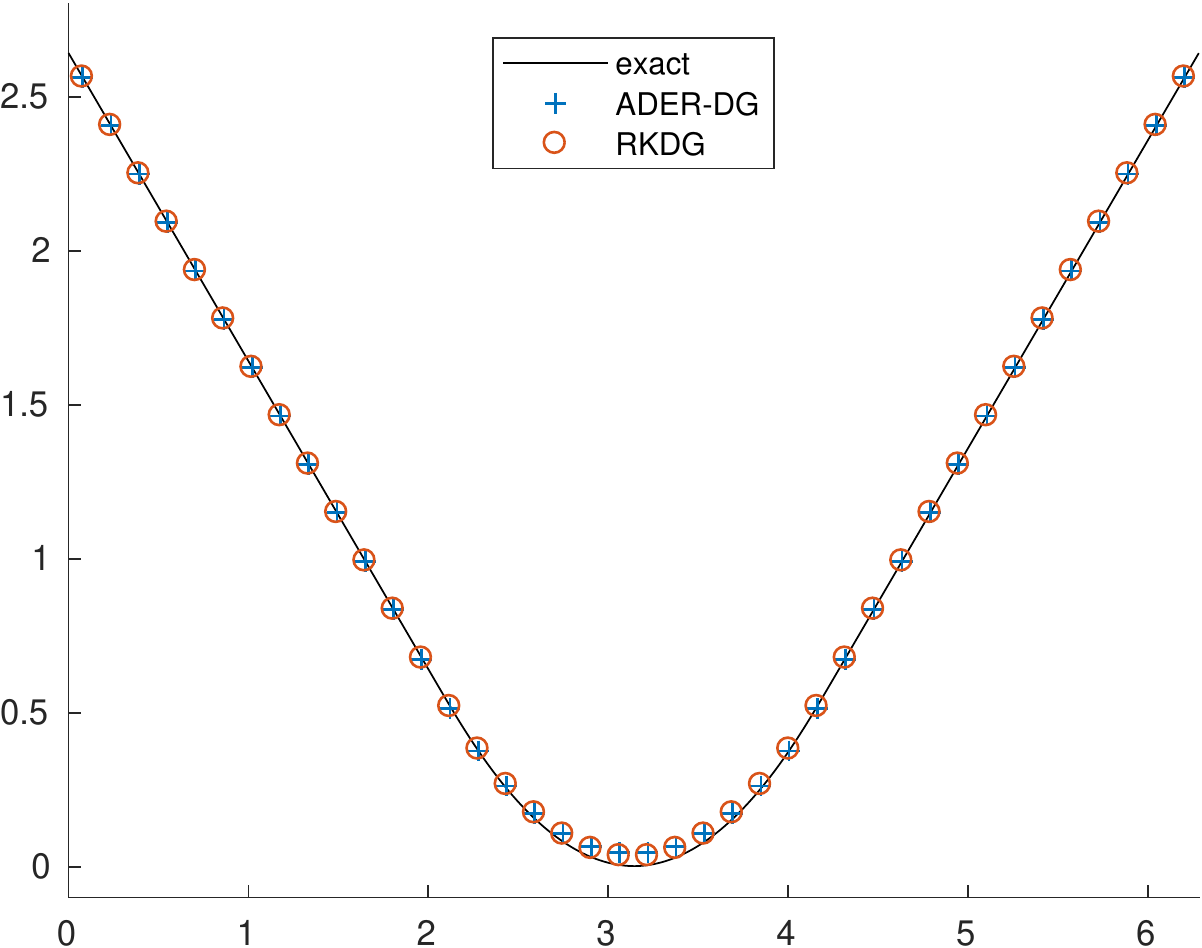}
  }
  \subfigure[$P^3,\text{CFL}=0.05$]{
    \includegraphics[width=0.45\textwidth]{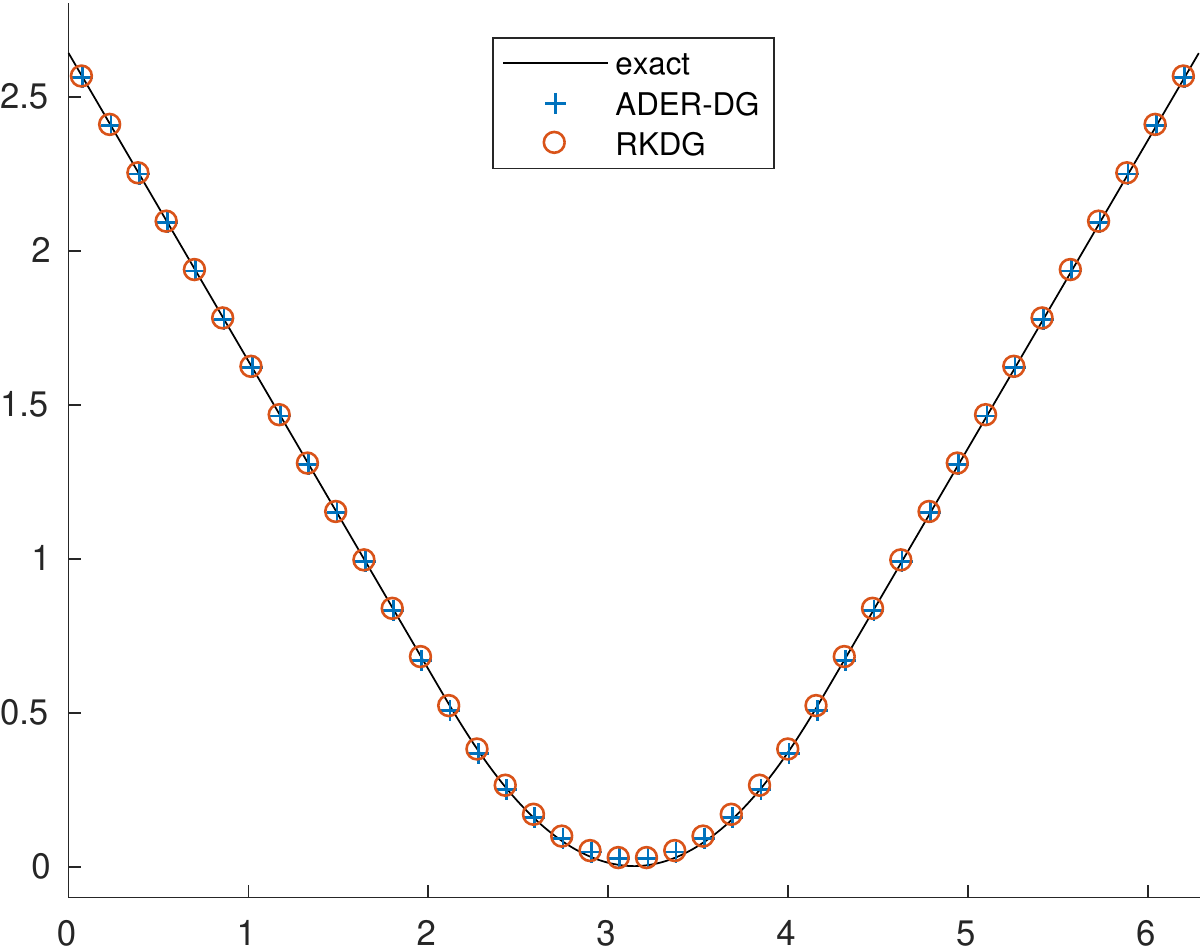}
  }
  \caption{Example \ref{ex:burgers_nonsmooth},~$N=40$. }
  \label{fig:burgers_nonsmooth}
\end{figure}

There is a rarefaction wave formed in the derivative of the exact solution,
thus the initial sharp corner at $x=\pi$ will be smeared out over time.
Fig. \ref{fig:burgers_nonsmooth} includes the $P^2$ and $P^3$ ADER-DG results with $N=40$.
Thanks to the penalty terms adding to the numerical fluxes, the results of the ADER-DG scheme
converge to the viscosity solution correctly.

\begin{example}
 {\rm \label{ex:-cos1D}
 	We solve one-dimensional HJ equation
 	\begin{equation*}
 	\varphi_t-\cos(\varphi_x+1)=0,\quad -1\leqslant x\leqslant 1,
 	\end{equation*}
with a nonconvex Hamiltonian, initial condition $\varphi(x,0)=-\cos(\pi x)$, and periodic boundary condition.
 }
\end{example}

We compute the solution up to $t=\dfrac{0.5}{\pi^2}$.
At this time, the solution is still smooth.
We list the errors and the orders of convergence in Table \ref{tab:cos1D}.
In the table, $(k+1)$-th order accuracy for $P^k$ polynomials can be observed.

\begin{table}[!ht]
  \centering
  \begin{tabular}{c|cccccc}
    \hline $N$ & $\ell^2$ error & Order & $\ell^1$ error & Order & $\ell^\infty$ error & Order \\
    \hline \multicolumn{7}{c}{$k=1,\text{CFL}=0.15$} \\ \hline
 10 & 1.082e-02 &  -   & 1.177e-02 &  -   & 2.009e-02 &  -   \\
 20 & 3.776e-03 & 1.52 & 3.199e-03 & 1.88 & 8.995e-03 & 1.16 \\
 40 & 7.522e-04 & 2.33 & 6.459e-04 & 2.31 & 2.063e-03 & 2.12 \\
 80 & 1.449e-04 & 2.38 & 1.380e-04 & 2.23 & 4.803e-04 & 2.10 \\
160 & 3.011e-05 & 2.27 & 2.918e-05 & 2.24 & 1.138e-04 & 2.08 \\
320 & 6.978e-06 & 2.11 & 7.261e-06 & 2.01 & 2.174e-05 & 2.39 \\
    \hline \multicolumn{7}{c}{$k=2,\text{CFL}=0.10$} \\ \hline
 10 & 1.345e-03 &  -   & 1.526e-03 &  -   & 2.355e-03 &  -   \\
 20 & 2.141e-04 & 2.65 & 2.359e-04 & 2.69 & 4.445e-04 & 2.41 \\
 40 & 2.968e-05 & 2.85 & 2.807e-05 & 3.07 & 7.981e-05 & 2.48 \\
 80 & 4.135e-06 & 2.84 & 3.728e-06 & 2.91 & 1.283e-05 & 2.64 \\
160 & 5.532e-07 & 2.90 & 4.891e-07 & 2.93 & 1.704e-06 & 2.91 \\
320 & 7.138e-08 & 2.95 & 6.098e-08 & 3.00 & 2.320e-07 & 2.88 \\
    \hline \multicolumn{7}{c}{$k=3,\text{CFL}=0.05$} \\ \hline
 10 & 3.630e-04 &  -   & 3.227e-04 &  -   & 1.099e-03 &  -   \\
 20 & 1.756e-05 & 4.37 & 1.551e-05 & 4.38 & 5.326e-05 & 4.37 \\
 40 & 1.436e-06 & 3.61 & 1.247e-06 & 3.64 & 4.616e-06 & 3.53 \\
 80 & 1.076e-07 & 3.74 & 8.771e-08 & 3.83 & 4.173e-07 & 3.47 \\
160 & 7.396e-09 & 3.86 & 5.831e-09 & 3.91 & 4.216e-08 & 3.31 \\
320 & 4.960e-10 & 3.90 & 3.817e-10 & 3.93 & 3.411e-09 & 3.63 \\
    \hline
  \end{tabular}
  \caption{Errors and orders of convergence for Example
  \ref{ex:-cos1D},~$t=\dfrac{0.5}{\pi^2}$}
  \label{tab:cos1D}
\end{table}

\begin{figure}[!ht]
  \centering
  \subfigure[$P^2$]{
    \includegraphics[width=0.45\textwidth]{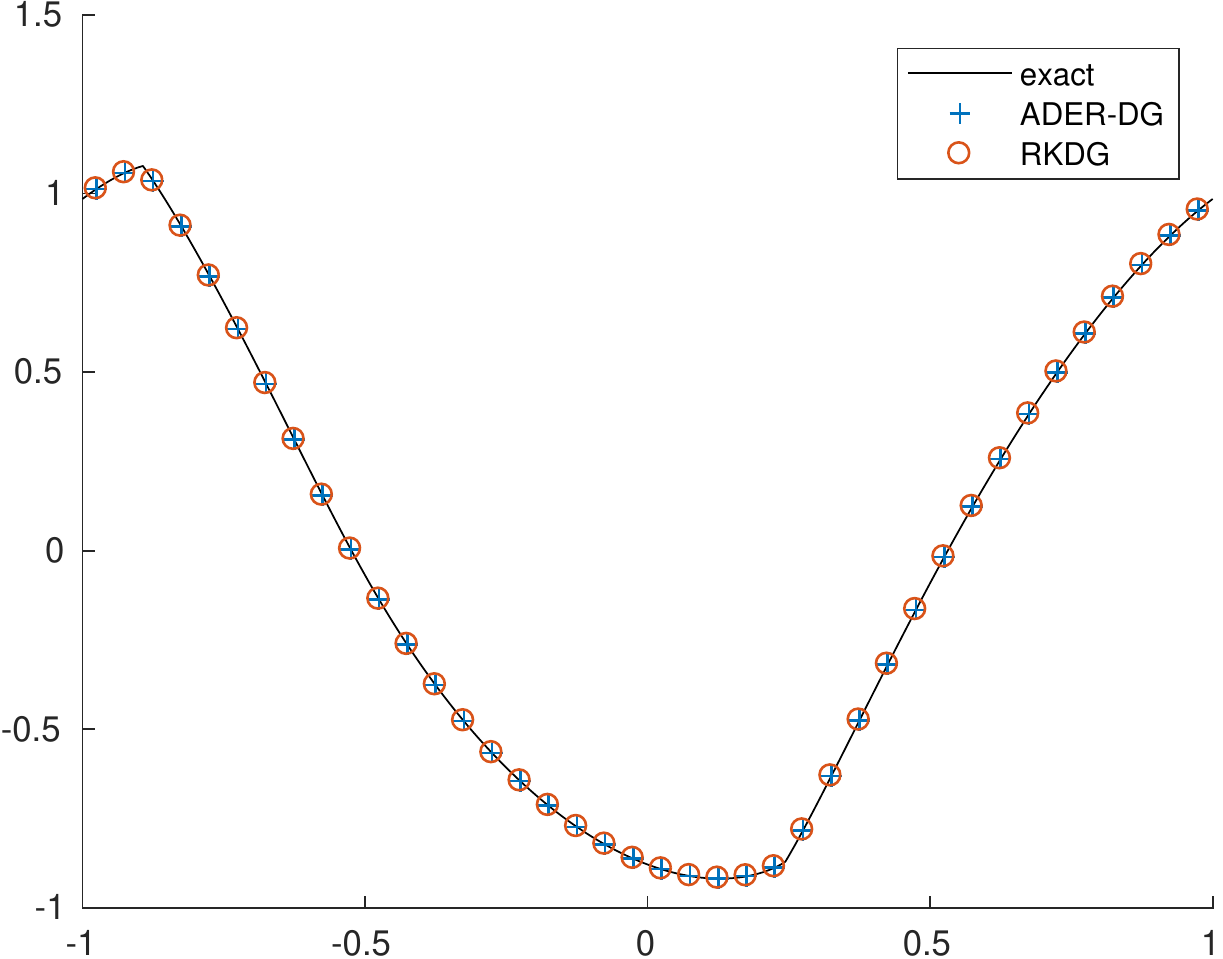}
  }
  \subfigure[$P^3$]{
    \includegraphics[width=0.45\textwidth]{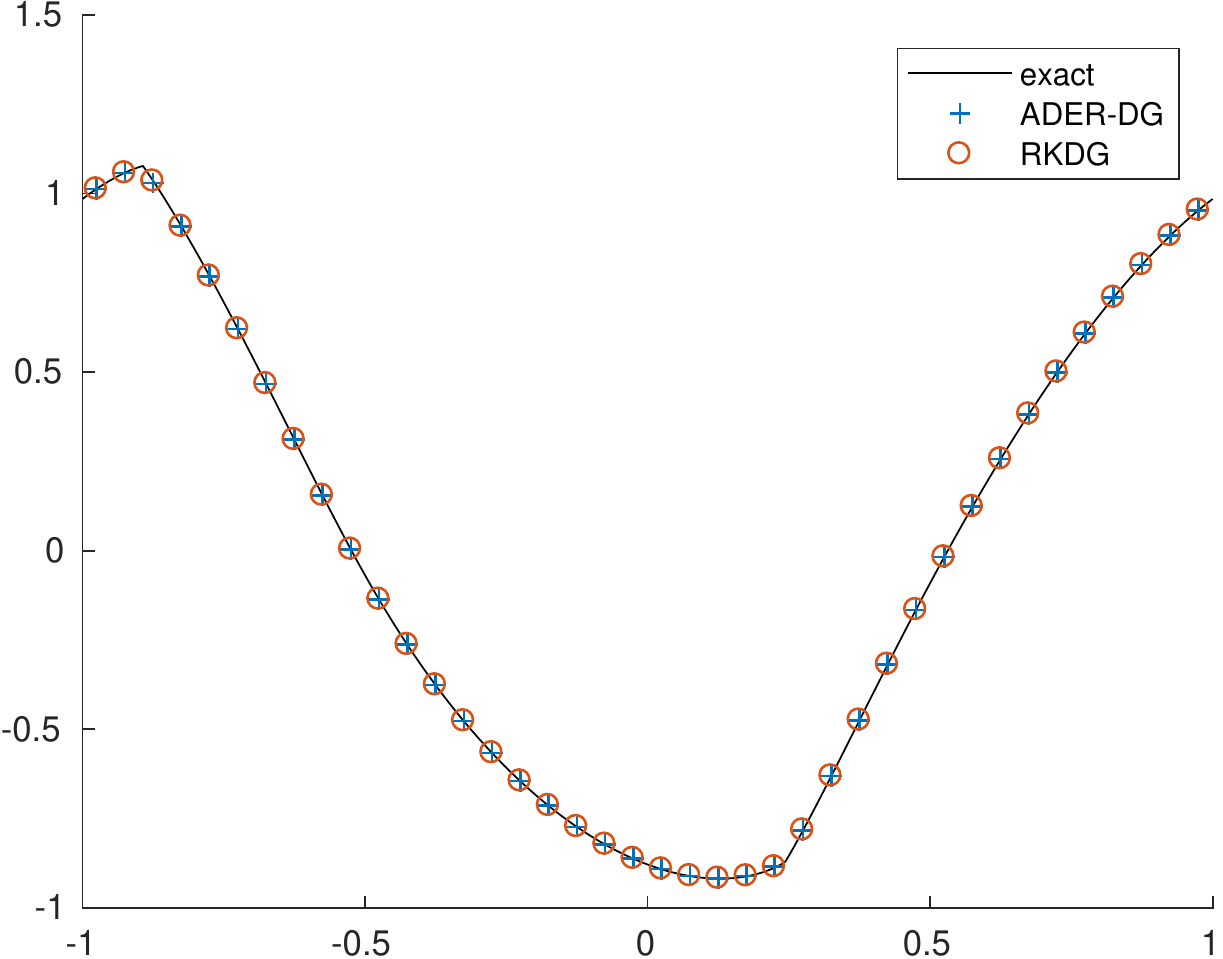}
  }
  \caption{Example \ref{ex:-cos1D},~$N=40$. }
  \label{fig:cos1D}
\end{figure}

Then we compute the solution up to $t=\dfrac{1.5}{\pi^2}$.
We plot the results of the ADER-DG scheme in Fig. \ref{fig:cos1D}.
In the figures, the kinks in the solution are clearly resolved by our scheme.

\begin{example}
{\rm \label{ex:nonconvex1D}
We solve one-dimensional Riemann problem
\begin{equation*}
\varphi_t+\dfrac{(\varphi_x^2-1)(\varphi_x^2-4)}{4}=0,\quad -1\leqslant x\leqslant 1,
\end{equation*}
with a nonconvex Hamiltonian, initial condition $\varphi(x,0)=-2\abs{x}$.
}
\end{example}

\begin{figure}[!ht]
  \centering
  \subfigure[$N=80$]{
    \includegraphics[width=0.45\textwidth]{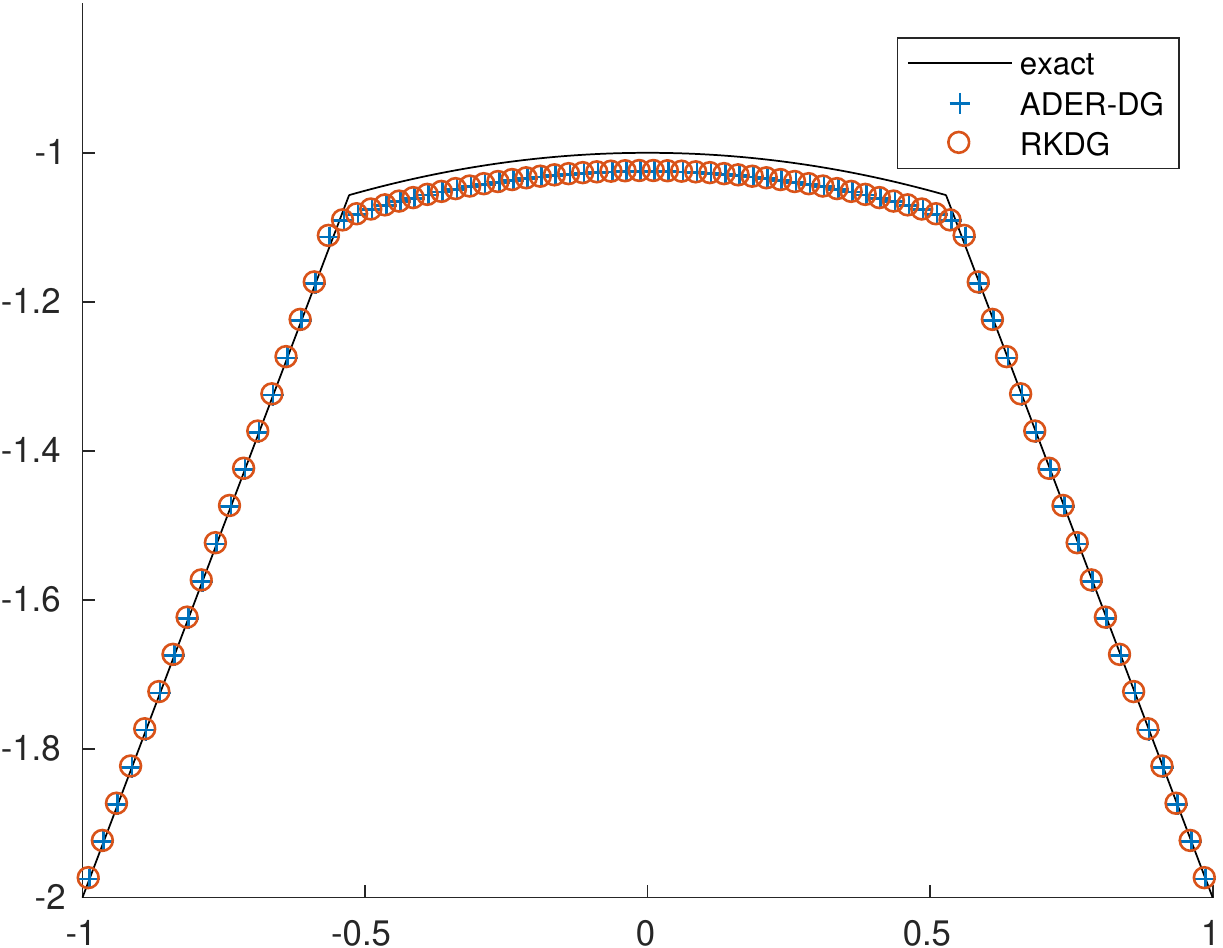}
  }
  \subfigure[$N=81$]{
    \includegraphics[width=0.45\textwidth]{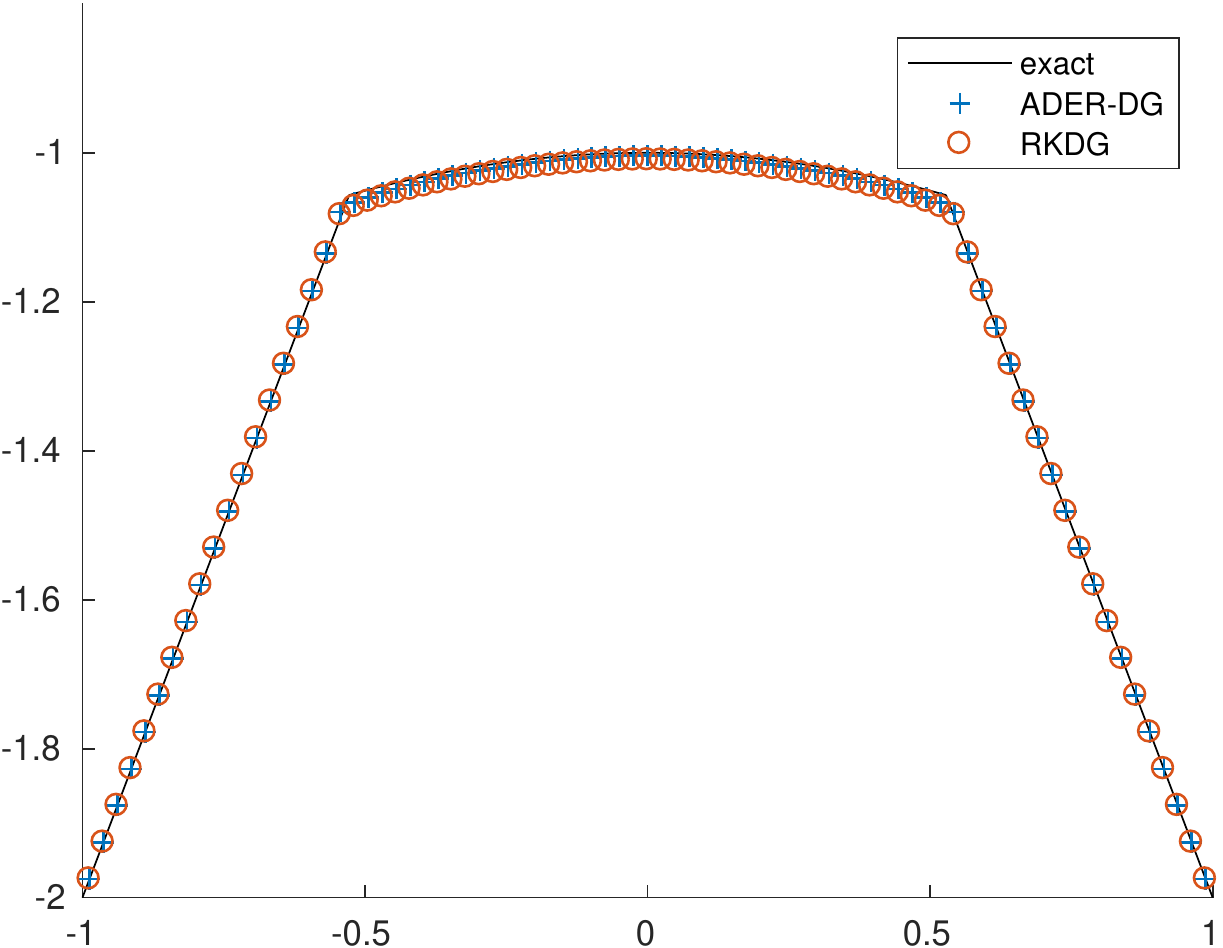}
  }
  \caption{Example \ref{ex:nonconvex1D},~$P^2,\text{CFL}=0.10$. Left: $N=80$, right: $N=81$. }
  \label{fig:nonconvex1D}
\end{figure}

The results at $t=1$ of the ADER-DG scheme with $N=80$ and $N=81$ are plotted in Fig. \ref{fig:nonconvex1D}.
It is a benchmark problem to test a numerical scheme's capability to capture the viscosity solution.
Similar as the RKDG scheme in \cite{Cheng2014}, Minmod limiter is used for the convergence to the entropy solution,
and the result with odd $N$ gives smaller errors.

\subsection{Two-dimensional results}
\begin{example}
{\rm \label{ex:linear12D}
	We solve the following linear problem with smooth variable coefficient \cite{Cheng2014}
	\begin{equation*}
	\varphi_t-y\varphi_x+x\varphi_y=0,\quad (x,y)\in [-1,1]^2,
	\end{equation*}
	with the initial condition
	\begin{equation*}
	\varphi(x,y,0)=\exp\left(-\dfrac{(x-0.4)^2+(y-0.4)^2}{2\sigma^2}\right),
	\end{equation*}
	and periodic boundary condition.
	The parameter $\sigma$ is $0.05$, and the computational time is $t=1$.
}
\end{example}

This problem describes a smooth solid body rotating around the origin.
We list the errors and the orders of convergence in Table \ref{tab:linear12D}.

\begin{table}[!ht]
  \centering
  \begin{tabular}{c|cccccc}
    \hline $N$ & $\ell^2$ error & Order & $\ell^1$ error & Order & $\ell^\infty$ error & Order \\
    \hline \multicolumn{7}{c}{$k=1,\text{CFL}=0.15$} \\ \hline
 10 & 6.631e-02 &  -   & 2.646e-02 &  -   & 5.453e-01 &  -   \\
 20 & 4.472e-02 & 0.57 & 1.296e-02 & 1.03 & 4.730e-01 & 0.21 \\
 40 & 2.085e-02 & 1.10 & 4.872e-03 & 1.41 & 2.831e-01 & 0.74 \\
 80 & 5.349e-03 & 1.96 & 1.066e-03 & 2.19 & 9.463e-02 & 1.58 \\
160 & 9.966e-04 & 2.42 & 1.877e-04 & 2.51 & 1.988e-02 & 2.25 \\
    \hline \multicolumn{7}{c}{$k=2,\text{CFL}=0.05$} \\ \hline
 10 & 4.743e-02 &  -   & 1.828e-02 &  -   & 4.470e-01 &  -   \\
 20 & 2.097e-02 & 1.18 & 5.772e-03 & 1.66 & 2.832e-01 & 0.66 \\
 40 & 3.030e-03 & 2.79 & 6.260e-04 & 3.20 & 5.288e-02 & 2.42 \\
 80 & 2.383e-04 & 3.67 & 4.720e-05 & 3.73 & 4.951e-03 & 3.42 \\
160 & 2.338e-05 & 3.35 & 4.744e-06 & 3.31 & 4.190e-04 & 3.56 \\
\hline \multicolumn{7}{c}{$k=3,\text{CFL}=0.05$} \\ \hline
 10 & 3.508e-02 &  -   & 1.222e-02 &  -   & 3.618e-01 &  -   \\
 20 & 6.505e-03 & 2.43 & 1.706e-03 & 2.84 & 7.337e-02 & 2.30 \\
 40 & 3.530e-04 & 4.20 & 7.418e-05 & 4.52 & 6.290e-03 & 3.54 \\
 80 & 1.682e-05 & 4.39 & 3.501e-06 & 4.41 & 5.178e-04 & 3.60 \\
160 & 9.300e-07 & 4.18 & 1.958e-07 & 4.16 & 2.627e-05 & 4.30 \\
    \hline
  \end{tabular}
  \caption{Errors and orders of convergence for Example
  \ref{ex:linear12D},~$t=1$}
  \label{tab:linear12D}
\end{table}

\begin{example}
{\rm \label{ex:linear22D}
  We solve the same problem as Example \ref{ex:linear12D},
  but with a unsmooth initial condition
  \begin{equation}
    \varphi(x,y,0)=\begin{cases}
      0, \quad & 0.3\leqslant r, \\
      0.3-r, \quad & 0.1<r<0.3, \\
      0.2, \quad & r\leqslant 0.1, \\
    \end{cases}
  \end{equation}
  where $r=\sqrt{(x-0.4)^2+(y-0.4)^2}$.
}
\end{example}

\begin{figure}[!ht]
  \centering
  \subfigure[$P^1$]{
    \includegraphics[width=0.3\textwidth]{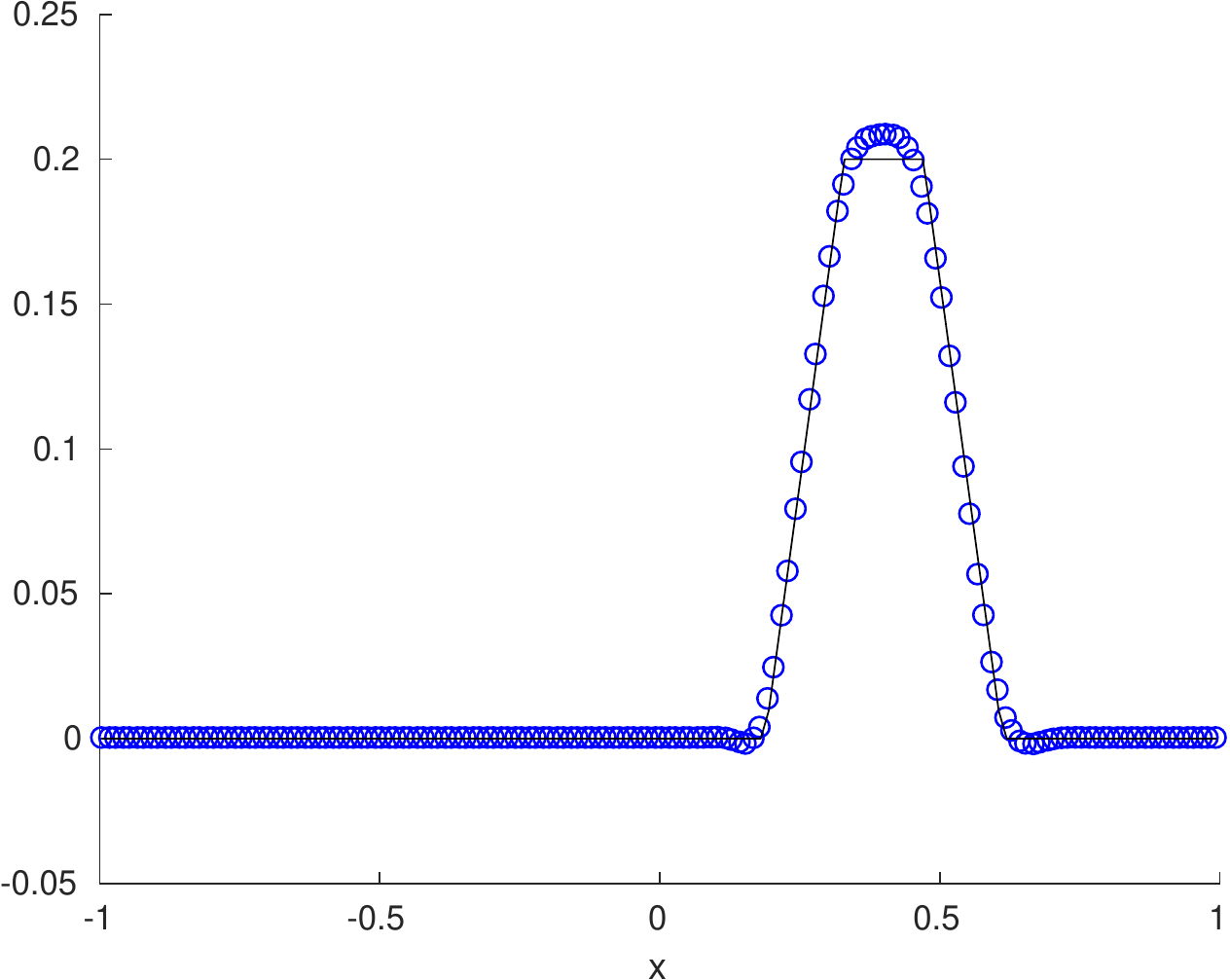}
  }
  \subfigure[$P^2$]{
    \includegraphics[width=0.3\textwidth]{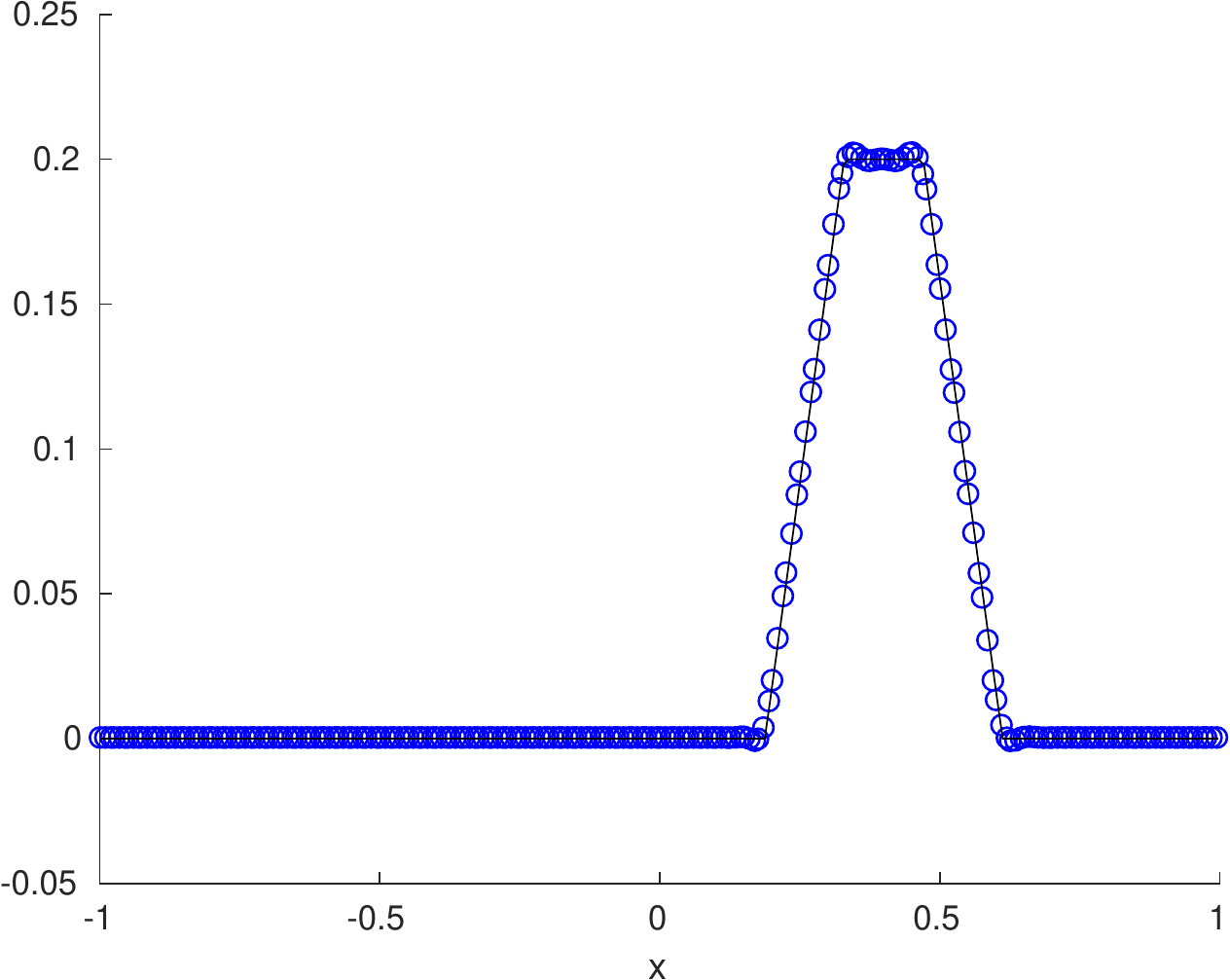}
  }
  \subfigure[$P^3$]{
    \includegraphics[width=0.3\textwidth]{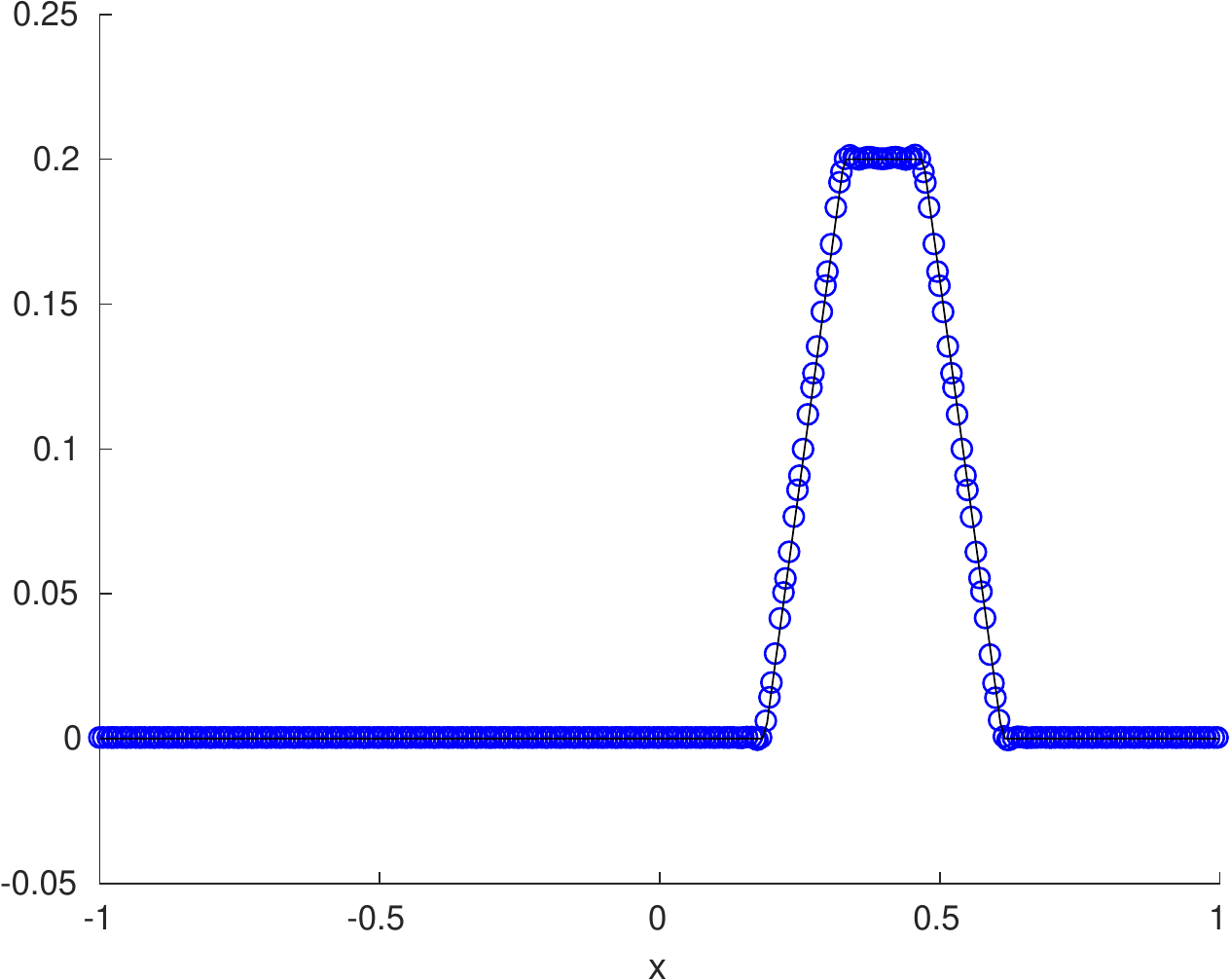}
  }
  \caption{The comparisons of $\varphi$ cut along the line $y=x$ for Example
  \ref{ex:linear22D},~$N=80$. Solid line is the exact solution and the
circles are numerical solutions obtained by the ADER-DG scheme with $P^1,P^2,P^3$ polynomials.}
  \label{fig:linear22D}
\end{figure}

The numerical results at $t=2\pi$ are provided in Table \ref{tab:linear22D}.
From the table, we can observe that the ADER-DG scheme is nearly first order,
because the initial condition is unsmooth.
But we can see from Fig. \ref{fig:linear22D} that, high order scheme can obtain better results.

\begin{table}[!ht]
  \centering
  \begin{tabular}{c|cccccc}
    \hline $N$ & $\ell^2$ error & Order & $\ell^1$ error & Order & $\ell^\infty$ error & Order \\
    \hline \multicolumn{7}{c}{$k=1,\text{CFL}=0.15$} \\ \hline
 10 & 3.354e-02 &  -   & 2.623e-02 &  -   & 1.176e-01 &  -   \\
 20 & 1.630e-02 & 1.04 & 1.198e-02 & 1.13 & 5.408e-02 & 1.12 \\
 40 & 6.029e-03 & 1.43 & 3.973e-03 & 1.59 & 2.052e-02 & 1.40 \\
 80 & 2.705e-03 & 1.16 & 1.512e-03 & 1.39 & 1.225e-02 & 0.74 \\
160 & 1.268e-03 & 1.09 & 5.870e-04 & 1.36 & 7.223e-03 & 0.76 \\
    \hline \multicolumn{7}{c}{$k=2,\text{CFL}=0.05$} \\ \hline
 10 & 1.336e-02 &  -   & 1.082e-02 &  -   & 3.939e-02 &  -   \\
 20 & 4.216e-03 & 1.66 & 2.883e-03 & 1.91 & 1.953e-02 & 1.01 \\
 40 & 1.964e-03 & 1.10 & 1.149e-03 & 1.33 & 9.025e-03 & 1.11 \\
 80 & 8.123e-04 & 1.27 & 3.794e-04 & 1.60 & 5.138e-03 & 0.81 \\
160 & 3.462e-04 & 1.23 & 1.243e-04 & 1.61 & 2.919e-03 & 0.82 \\
    \hline \multicolumn{7}{c}{$k=3,\text{CFL}=0.05$} \\ \hline
 10 & 5.797e-03 &  -   & 4.822e-03 &  -   & 2.688e-02 &  -   \\
 20 & 2.525e-03 & 1.20 & 1.656e-03 & 1.54 & 1.103e-02 & 1.29 \\
 40 & 9.705e-04 & 1.38 & 5.182e-04 & 1.68 & 5.527e-03 & 1.00 \\
 80 & 3.982e-04 & 1.29 & 1.643e-04 & 1.66 & 3.026e-03 & 0.87 \\
160 & 1.639e-04 & 1.28 & 5.113e-05 & 1.68 & 1.721e-03 & 0.81 \\
    \hline
  \end{tabular}
  \caption{Errors and orders of convergence for Example
  \ref{ex:linear22D},~$t=2\pi$}
  \label{tab:linear22D}
\end{table}

\begin{example}
{\rm \label{ex:burgers2D}
	We solve two-dimensional Burgers' equation
	\begin{equation*}
	\varphi_t+\dfrac{(\varphi_x+\varphi_y+1)^2}{2}=0,\quad (x,y)\in [-2,2]^2,
	\end{equation*}
with a smooth initial condition  $\varphi(x,y,0)=-\cos(\frac{\pi}{2}(x+y))$, and periodic boundary condition.
}
\end{example}

\begin{figure}[!ht]
  \centering
  \subfigure[$t=\dfrac{0.5}{\pi^2}$]{
    \includegraphics[width=0.45\textwidth]{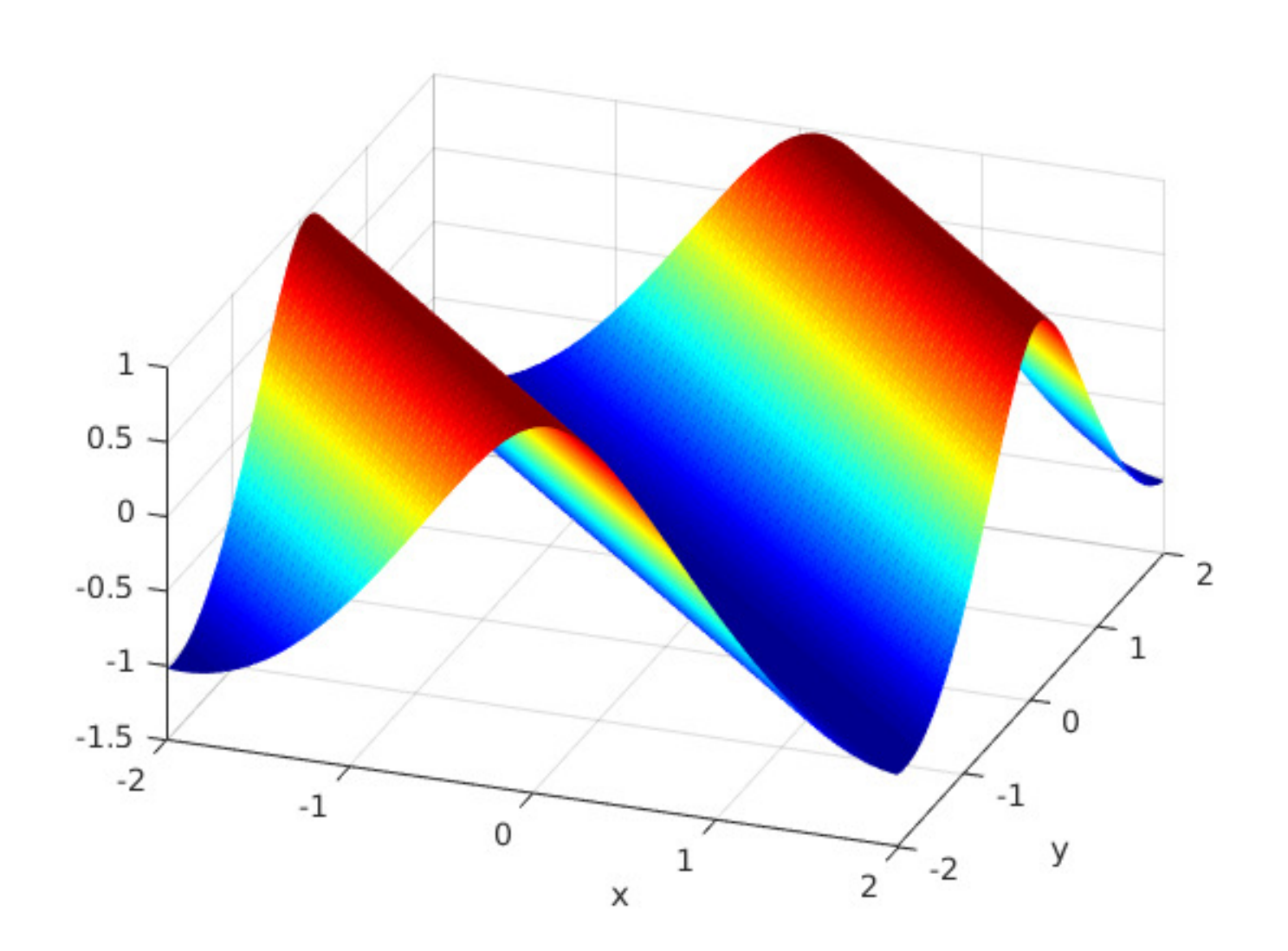}
  }
  \subfigure[$t=\dfrac{1.5}{\pi^2}$]{
    \includegraphics[width=0.45\textwidth]{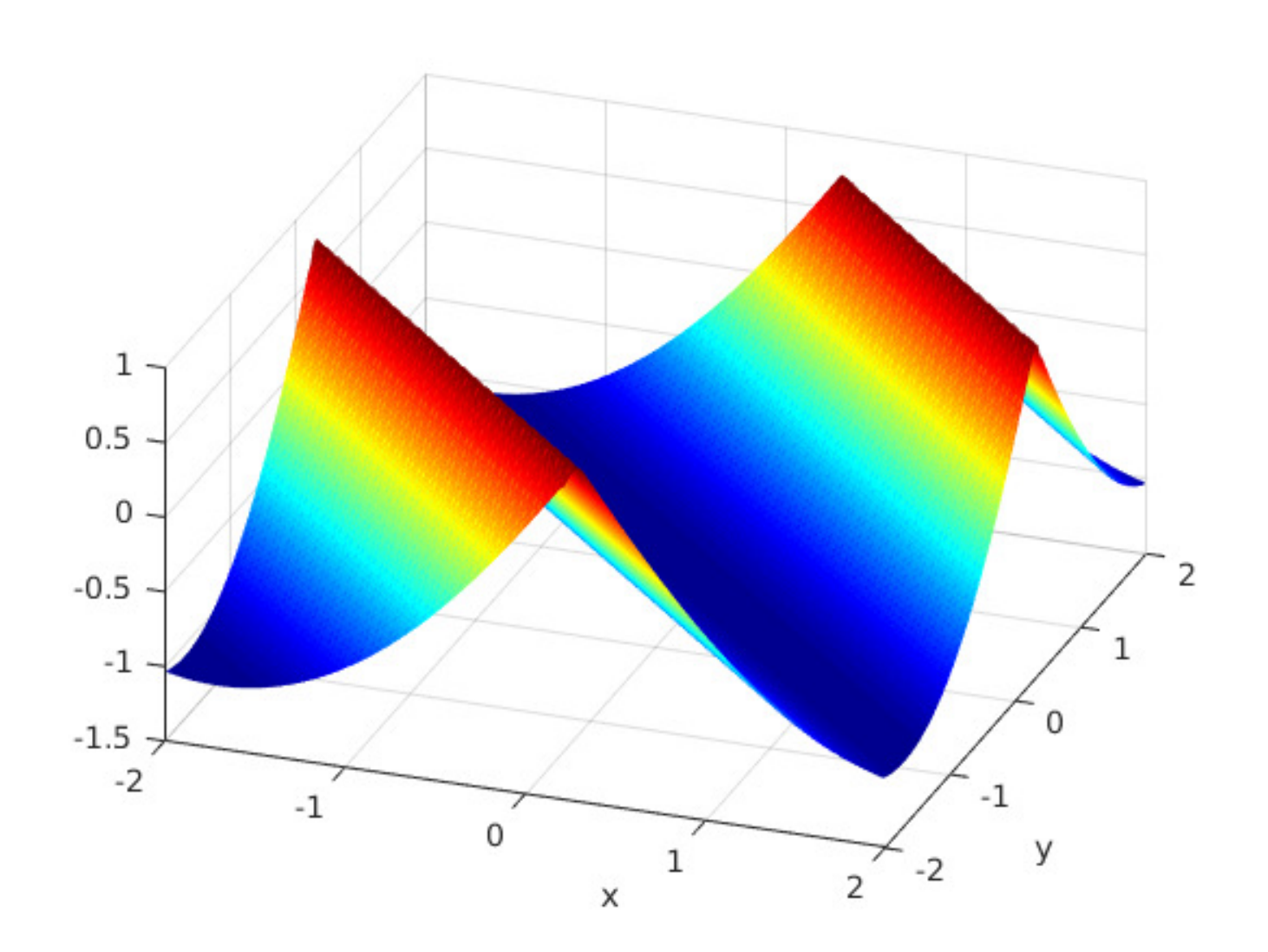}
  }
  \caption{Example \ref{ex:burgers2D},~$P^2,N=40$}
  \label{fig:burgers2D}
\end{figure}

We compute  the solution until $t=\dfrac{0.5}{\pi^2}$. It is smooth at this time.
We give the numerical errors and the orders of convergence in Table \ref{tab:burgers2D}.
It is clearly that our scheme can achieve $(k+1)$-th order of convergence for $P^k$ polynomial.
We also compute the same equation until $t=\dfrac{1.5}{\pi^2}$,
and the discontinuous derivative has already appeared in the solution.
We plot the results in Fig. \ref{fig:burgers2D},
from which we can observe good resolutions of the ADER-DG scheme for this example.

\begin{table}[!ht]
  \centering
  \begin{tabular}{c|cccccc}
    \hline $N$ & $\ell^2$ error & Order & $\ell^1$ error & Order & $\ell^\infty$ error & Order \\
    \hline \multicolumn{7}{c}{$k=1,\text{CFL}=0.15$} \\ \hline
 10 & 1.121e-01 &  -   & 3.599e-01 &  -   & 9.131e-02 &  -   \\
 20 & 2.855e-02 & 1.97 & 8.840e-02 & 2.03 & 2.609e-02 & 1.81 \\
 40 & 7.174e-03 & 1.99 & 2.200e-02 & 2.01 & 6.557e-03 & 1.99 \\
 80 & 1.808e-03 & 1.99 & 5.467e-03 & 2.01 & 1.800e-03 & 1.86 \\
160 & 4.579e-04 & 1.98 & 1.361e-03 & 2.01 & 4.963e-04 & 1.86 \\
320 & 1.160e-04 & 1.98 & 3.391e-04 & 2.00 & 1.326e-04 & 1.90 \\
    \hline \multicolumn{7}{c}{$k=2,\text{CFL}=0.10$} \\ \hline
 10 & 1.970e-02 &  -   & 5.484e-02 &  -   & 2.837e-02 &  -   \\
 20 & 2.729e-03 & 2.85 & 7.248e-03 & 2.92 & 3.011e-03 & 3.24 \\
 40 & 3.568e-04 & 2.94 & 8.996e-04 & 3.01 & 4.757e-04 & 2.66 \\
 80 & 4.662e-05 & 2.94 & 1.157e-04 & 2.96 & 6.699e-05 & 2.83 \\
160 & 6.024e-06 & 2.95 & 1.474e-05 & 2.97 & 9.528e-06 & 2.81 \\
320 & 7.714e-07 & 2.97 & 1.869e-06 & 2.98 & 1.296e-06 & 2.88 \\
    \hline \multicolumn{7}{c}{$k=3,\text{CFL}=0.05$} \\ \hline
 10 & 5.118e-03 &  -   & 9.484e-03 &  -   & 1.798e-02 &  -   \\
 20 & 2.893e-04 & 4.14 & 5.677e-04 & 4.06 & 8.581e-04 & 4.39 \\
 40 & 1.932e-05 & 3.90 & 3.588e-05 & 3.98 & 6.199e-05 & 3.79 \\
 80 & 1.258e-06 & 3.94 & 2.250e-06 & 4.00 & 4.992e-06 & 3.63 \\
160 & 8.186e-08 & 3.94 & 1.420e-07 & 3.99 & 3.548e-07 & 3.81 \\
320 & 5.299e-09 & 3.95 & 9.022e-09 & 3.98 & 2.391e-08 & 3.89 \\
    \hline
  \end{tabular}
  \caption{Errors and orders of convergence for Example
  \ref{ex:burgers2D},~$t=\dfrac{0.5}{\pi^2}$}
  \label{tab:burgers2D}
\end{table}

\begin{example}
{\rm \label{ex:-cos2D}
	We solve the following two-dimensional equation with nonconvex Hamiltonian
	\begin{equation*}
	\varphi_t-\cos(\varphi_x+\varphi_y+1)=0,\quad (x,y)\in [-2,2]^2,
	\end{equation*}
	with initial condition $\varphi(x,y,0)=-\cos(\frac{\pi}{2}(x+y))$, and periodic boundary condition.
}
\end{example}

\begin{figure}[!ht]
  \centering
  \subfigure[$t=\dfrac{0.5}{\pi^2}$]{
    \includegraphics[width=0.45\textwidth]{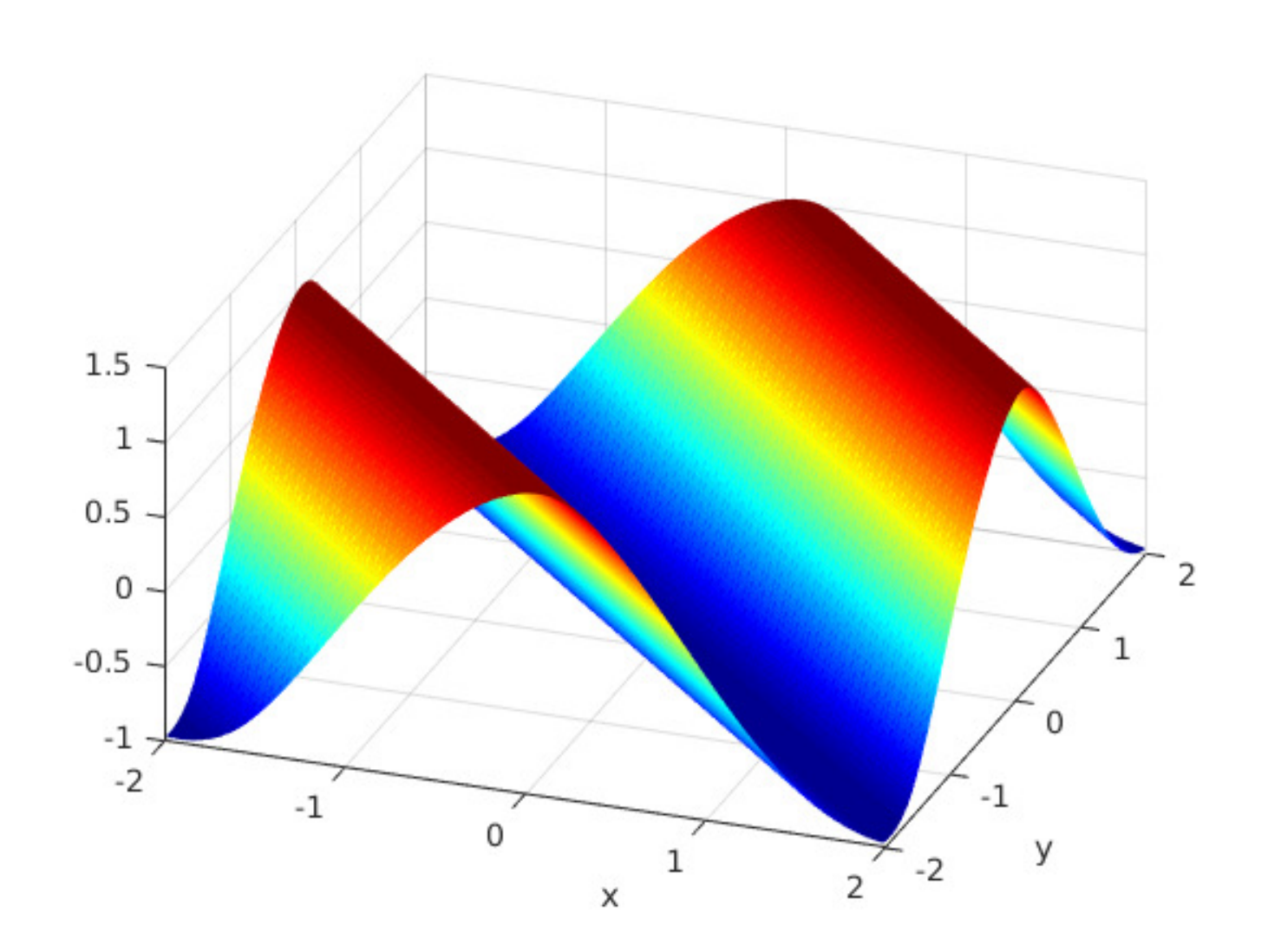}
  }
  \subfigure[$t=\dfrac{1.5}{\pi^2}$]{
    \includegraphics[width=0.45\textwidth]{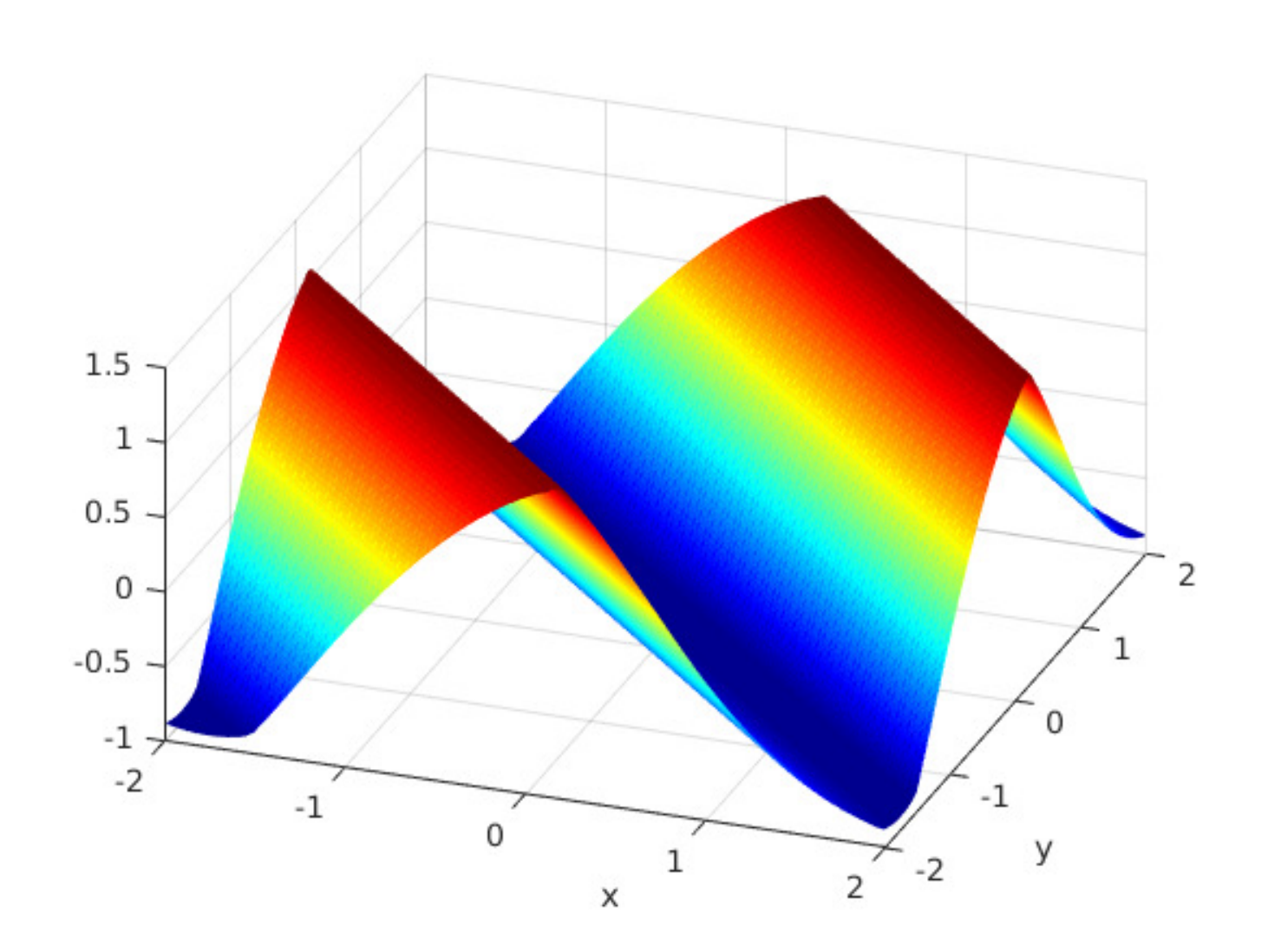}
  }
  \caption{Example \ref{ex:-cos2D},$P^2,N=40$}
  \label{fig:-cos2D}
\end{figure}

The solution is still smooth at $t=\dfrac{0.5}{\pi^2}$, and
the numerical errors and the orders of convergence at this time are listed in Table \ref{tab:-cos2D}.
We also compute the same equation until $t=\dfrac{1.5}{\pi^2}$,
and singular features develop in the solution.
The results of the ADER-DG scheme are shown in Fig. \ref{fig:-cos2D}.

\begin{table}[!ht]
  \centering
  \begin{tabular}{c|cccccc}
    \hline $N$ & $\ell^2$ error & Order & $\ell^1$ error & Order & $\ell^\infty$ error & Order \\
    \hline \multicolumn{7}{c}{$k=1,\text{CFL}=0.15$ } \\ \hline
 10 & 1.011e-01 &  -   & 3.371e-01 &  -   & 8.475e-02 &  -   \\
 20 & 2.649e-02 & 1.93 & 8.548e-02 & 1.98 & 2.185e-02 & 1.96 \\
 40 & 6.610e-03 & 2.00 & 2.144e-02 & 2.00 & 6.102e-03 & 1.84 \\
 80 & 1.648e-03 & 2.00 & 5.343e-03 & 2.00 & 1.630e-03 & 1.90 \\
160 & 4.156e-04 & 1.99 & 1.335e-03 & 2.00 & 4.330e-04 & 1.91 \\
320 & 1.030e-04 & 2.01 & 3.326e-04 & 2.00 & 1.127e-04 & 1.94 \\
    \hline \multicolumn{7}{c}{$k=2,\text{CFL}=0.10$ } \\ \hline
 10 & 3.503e-02 &  -   & 1.051e-01 &  -   & 2.916e-02 &  -   \\
 20 & 4.819e-03 & 2.86 & 1.456e-02 & 2.85 & 4.849e-03 & 2.59 \\
 40 & 5.946e-04 & 3.02 & 1.740e-03 & 3.07 & 6.406e-04 & 2.92 \\
 80 & 7.574e-05 & 2.97 & 2.102e-04 & 3.05 & 9.473e-05 & 2.76 \\
160 & 9.764e-06 & 2.96 & 2.629e-05 & 3.00 & 1.618e-05 & 2.55 \\
320 & 1.258e-06 & 2.96 & 3.303e-06 & 2.99 & 2.377e-06 & 2.77 \\
    \hline \multicolumn{7}{c}{$k=3,\text{CFL}=0.05$ } \\ \hline
 10 & 7.142e-03 &  -   & 1.732e-02 &  -   & 9.557e-03 &  -   \\
 20 & 8.431e-04 & 3.08 & 1.703e-03 & 3.35 & 1.887e-03 & 2.34 \\
 40 & 5.677e-05 & 3.89 & 1.004e-04 & 4.08 & 2.464e-04 & 2.94 \\
 80 & 3.687e-06 & 3.94 & 6.144e-06 & 4.03 & 2.008e-05 & 3.62 \\
160 & 2.400e-07 & 3.94 & 3.785e-07 & 4.02 & 1.323e-06 & 3.92 \\
320 & 1.541e-08 & 3.96 & 2.378e-08 & 3.99 & 8.863e-08 & 3.90 \\
    \hline
  \end{tabular}
  \caption{Errors and orders of convergence for Example
  \ref{ex:-cos2D},~$t=\dfrac{0.5}{\pi^2}$}
  \label{tab:-cos2D}
\end{table}

\begin{example}
{\rm \label{ex:controlling_optimal}
	We solve the following problem from optimal control
	\begin{equation*}
	\varphi_t+\sin(y)u+(\sin(x)+\text{sign}(v))v-\dfrac12\sin^2(y)+\cos(x)-1=0,\quad (x,y)\in [-\pi,\pi]^2,
	\end{equation*}
	with initial condition $\varphi(x,y,0)=0$, and periodic boundary condition.
}
\end{example}

\begin{figure}[!ht]
  \centering
  \subfigure[$\varphi$]{
    \includegraphics[width=0.45\textwidth]{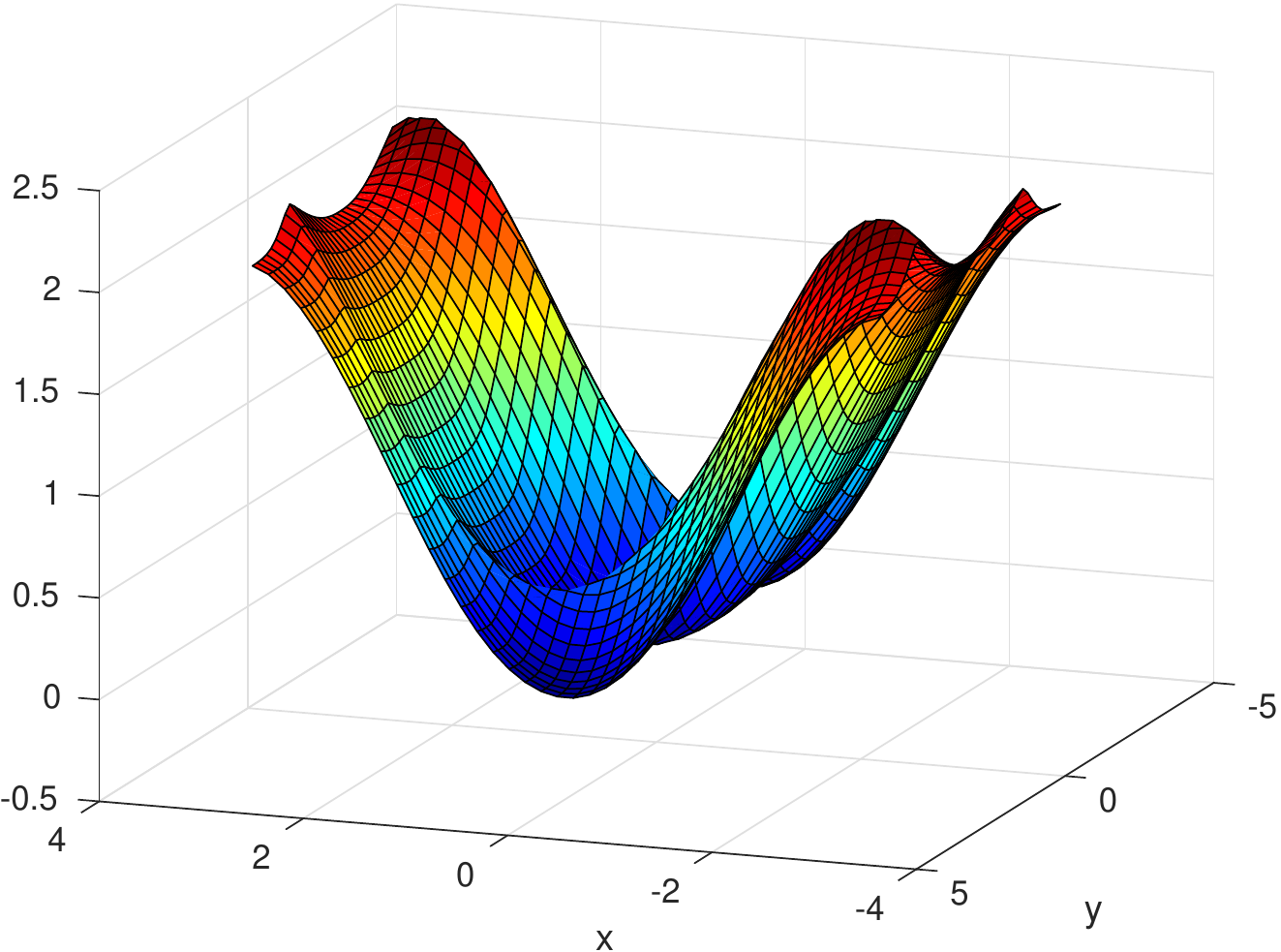}
  }
  \subfigure[$\varphi_y$]{
    \includegraphics[width=0.45\textwidth]{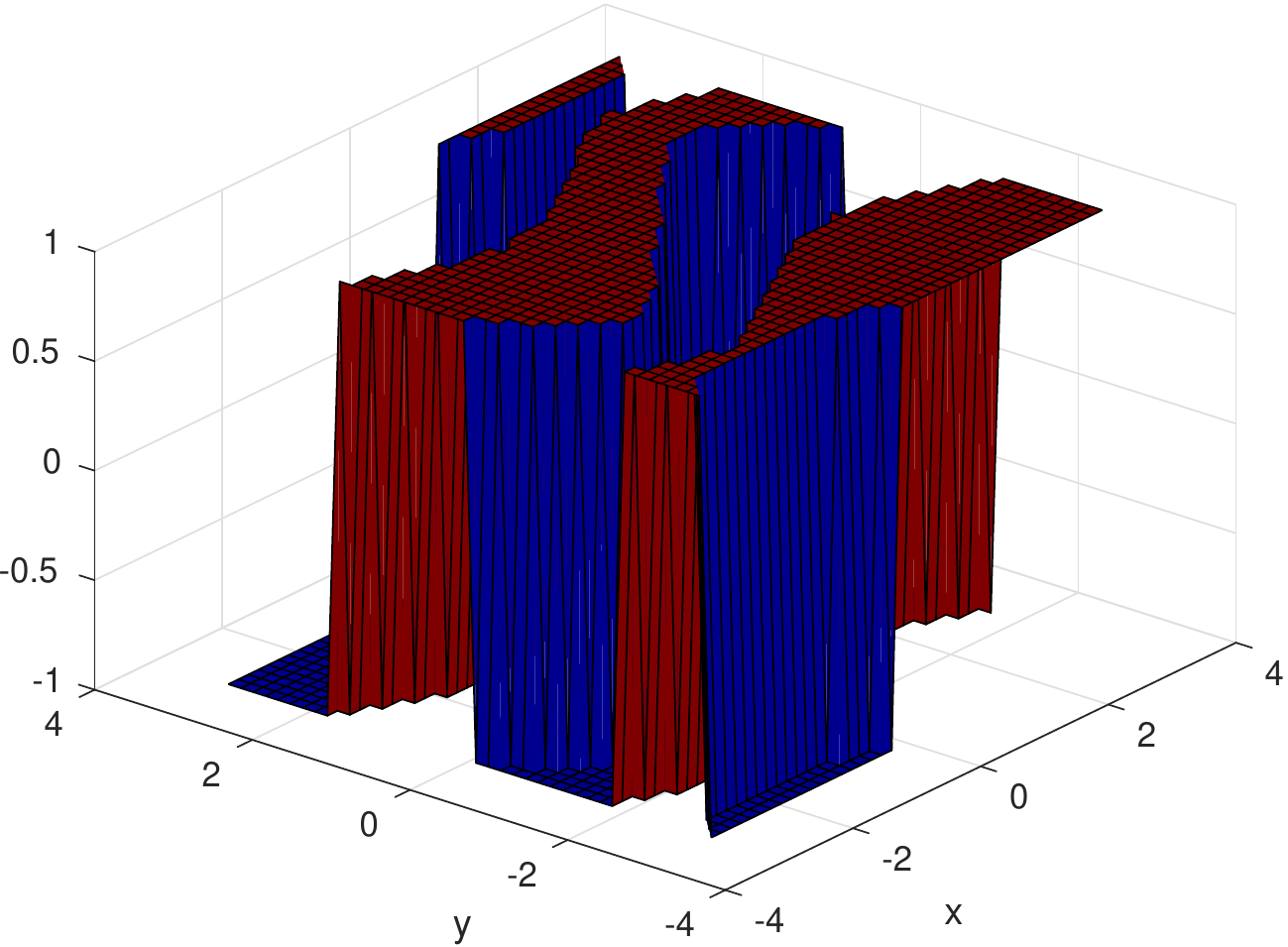}
  }
  \caption{Example \ref{ex:controlling_optimal},~$P^2,N=40,\text{CFL}=0.10$}
  \label{fig:controlling_optimal}
\end{figure}

The results obtained by the ADER-DG scheme with $P^2$ and $N=40$ are plotted in Fig. \ref{fig:controlling_optimal}.
We can see that our scheme can simulate the problem well.

\begin{example}
{\rm \label{ex:2DRiemann}
	We solve the following two-dimensional Riemann problem
	\begin{equation*}
	\varphi_t+\sin(\varphi_x+\varphi_y)=0,\quad (x,y)\in [-1,1]^2,
	\end{equation*}
	with initial condition $\varphi(x,y,0)=\pi(\abs{y}-\abs{x})$.
}
\end{example}

\begin{figure}[!ht]
  \centering
  \subfigure[$P^1,\text{CFL}=0.15$]{
    \includegraphics[width=0.45\textwidth]{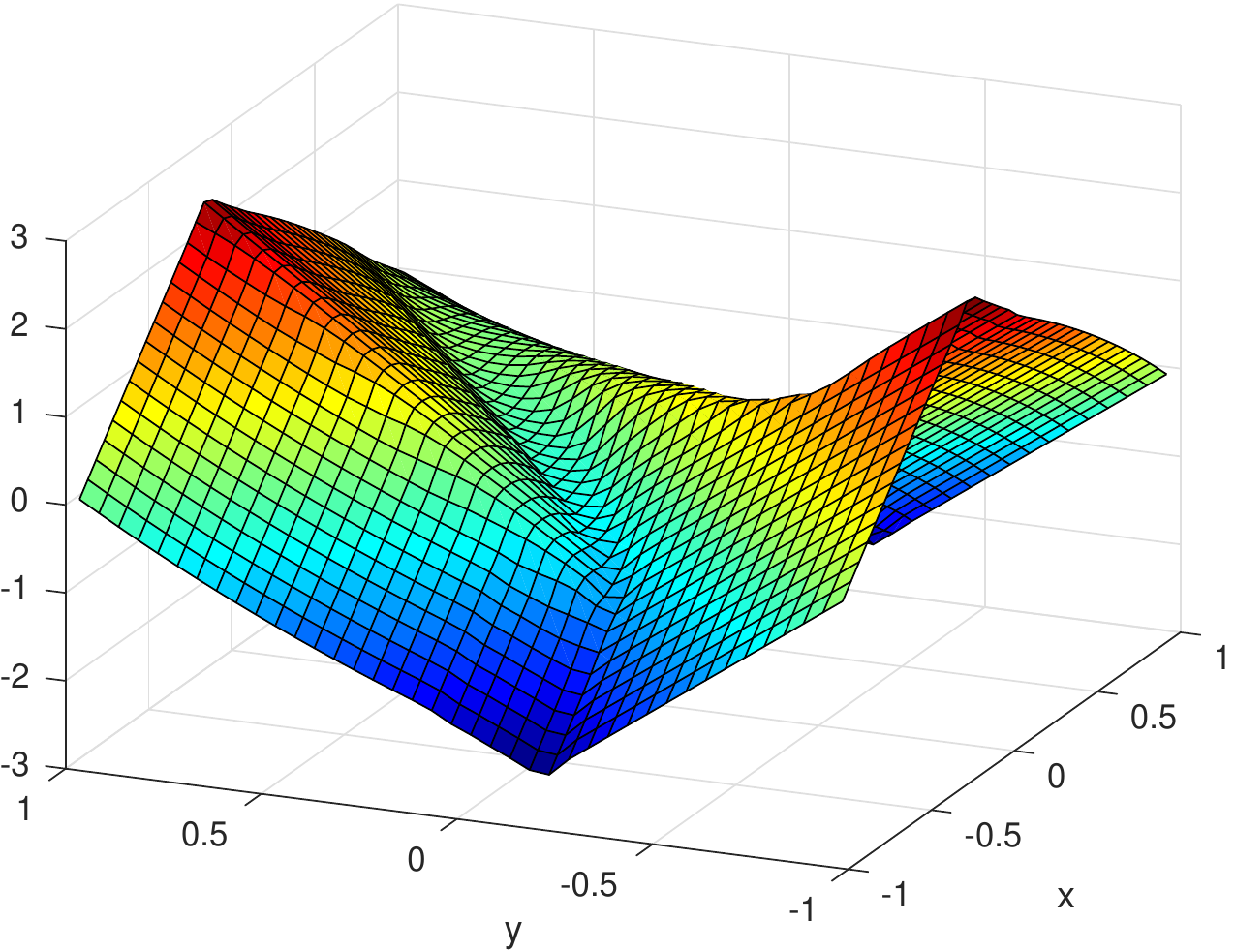}
  }
  \subfigure[$P^2,\text{CFL}=0.10$]{
    \includegraphics[width=0.45\textwidth]{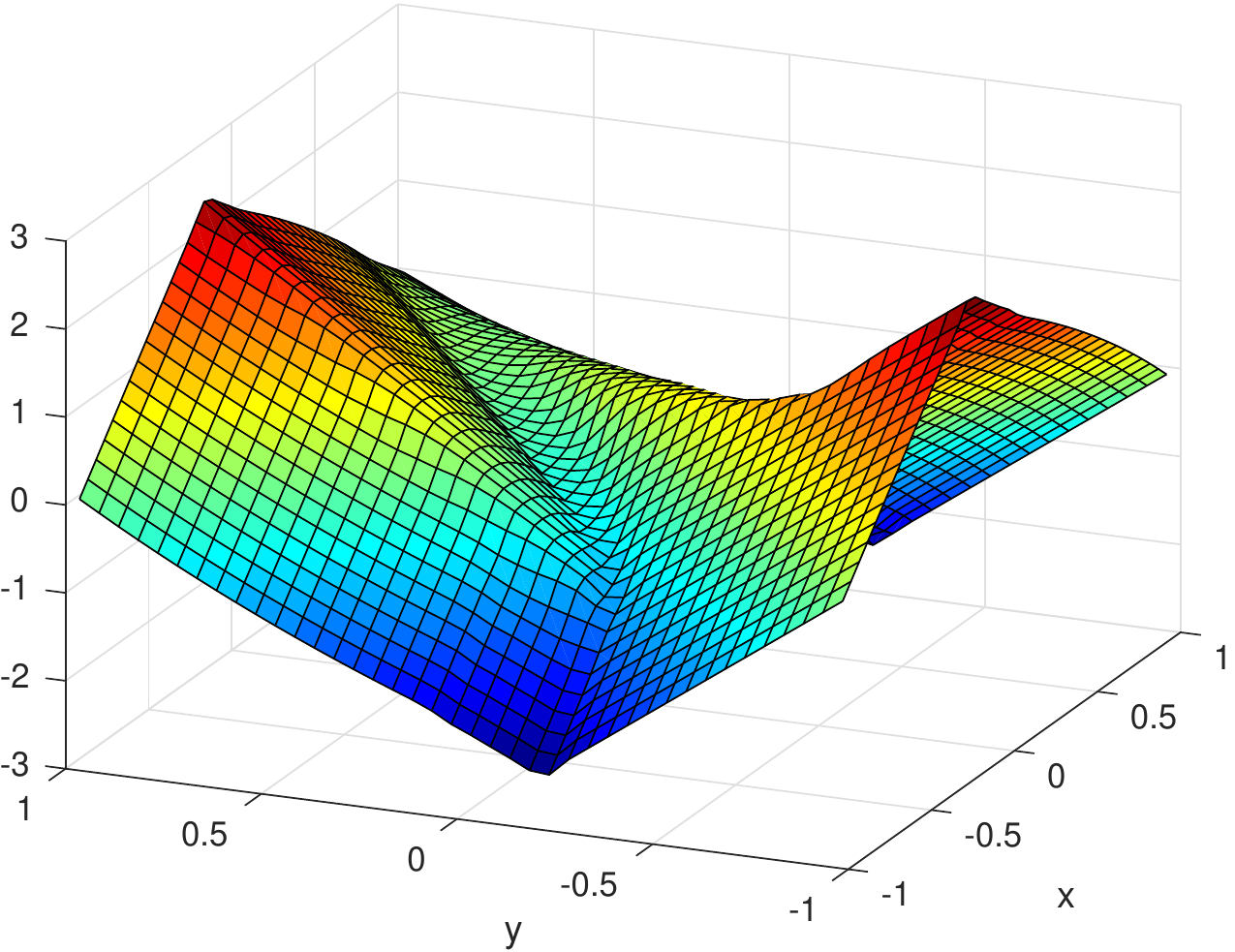}
  }
  \caption{Example \ref{ex:2DRiemann},~$N=40$}
  \label{fig:2DRiemann}
\end{figure}

We need limiters in this example to have its convergence to the viscosity solution.
The results of $P^1$ and $P^2$ ADER-DG scheme with $N=40$ are given in Fig. \ref{fig:2DRiemann}.
Our results are nearly the same as that in \cite{Cheng2014}.

\begin{example}
{\rm \label{ex:propagating_surface}
	We solve the problem of a propagating surface, which is a special case of the example in \cite{Osher1988}
	\begin{equation*}
	\varphi_t-\sqrt{\varphi_x^2+\varphi_y^2+1}=0,\quad (x,y)\in [0,1]^2,
	\end{equation*}
	with initial condition $\varphi(x,y,0)=1-\dfrac14(\cos(2\pi x)-1)(\cos(2\pi y)-1)$.
}
\end{example}

\begin{figure}[!ht]
  \centering
  \subfigure[$P^2,\text{CFL}=0.10$]{
    \includegraphics[width=0.45\textwidth]{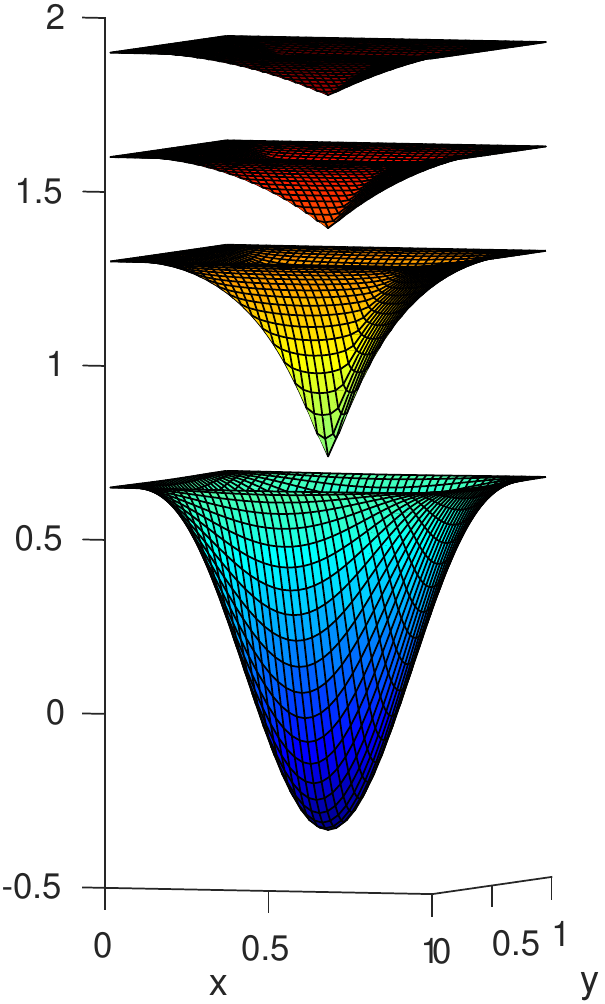}
  }
  \subfigure[$P^3,\text{CFL}=0.05$]{
    \includegraphics[width=0.45\textwidth]{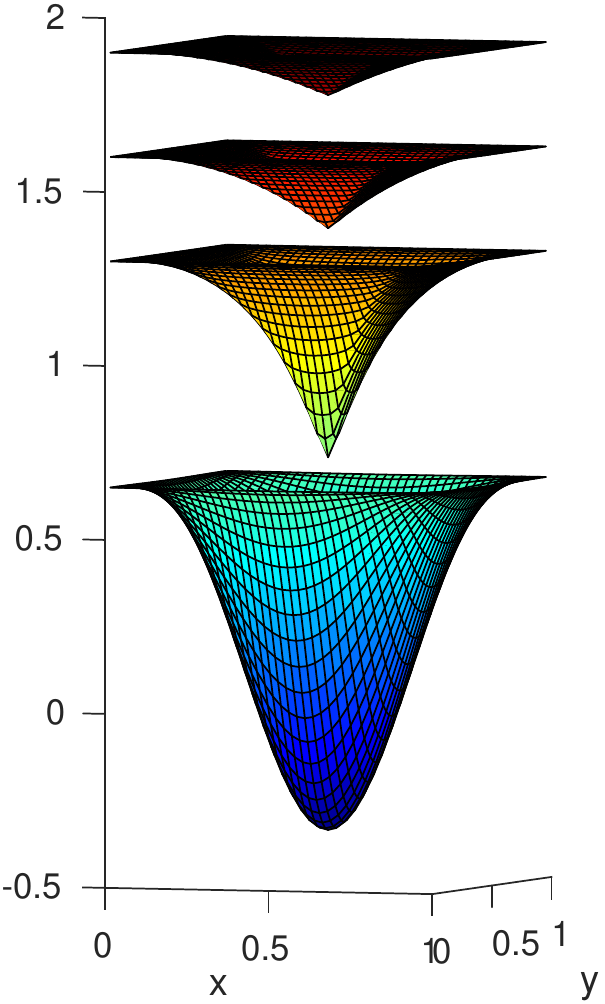}
  }
  \caption{Example \ref{ex:propagating_surface},~$N=41$}
  \label{fig:propagating_surface}
\end{figure}

We output the results at $t=0,0.3,0.6,0.9$, and plot them in Fig. \ref{fig:propagating_surface}.
The result at $t=0$ is shifted down to show the detail of the solution at later time.

At the end of this section,  presents a comparison of the CPU times of the  ADER-DG and RKDG schemes when they are applied to the above three two-dimensional examples, Examples \ref{ex:linear12D}, \ref{ex:burgers2D}, and \ref{ex:-cos2D}.
To make the fair comparison, we take the largest CFL number of both schemes,
although we have to use a slightly smaller CFL number for the ADER-DG scheme
in some cases.
The results given in Table \ref{tab:CPU} show that
the average CPU times of the ADER-DG scheme are about $25\%,21\%,14\%$ of the RKDG scheme
at second, third and fourth order, respectively, and
the ADER-DG scheme is more efficient than the RKDG scheme.

\begin{table}[!ht]
	\centering
	\begin{tabular}{c|r|rrr}\hline
		& & \tabincell{c}{Example \ref{ex:linear12D}\\ with $N=160$} &
		\tabincell{c}{Example \ref{ex:burgers2D}\\ with $N=320$} &
		\tabincell{c}{Example \ref{ex:-cos2D}\\ with $N=320$} \\ \hline
    \multirow{2}{*}{$P^1$} & ADER-DG & 3.95 & 3.43 & 1.32  \\
    & RKDG & 10.9 & 9.69 & 4.11 \\ \hline
    \multirow{2}{*}{$P^2$} & ADER-DG & 19.5 & 8.67 & 3.90  \\
    & RKDG & 80.2 & 47.9 & 26.8 \\ \hline
    \multirow{2}{*}{$P^3$} & ADER-DG & 45.3 & 36.3 & 15.7  \\
    & RKDG & 405 & 185 & 87.7 \\ \hline
	\end{tabular}
	\caption{The CPU times (second) of the ADER-DG scheme and the RKDG scheme. }
	\label{tab:CPU}
\end{table}

\section{Conclusion}\label{sec:conclusion}
An efficient ADER-DG scheme {was} presented to directly solve the Hamilton-Jacobi equations.
The ADER-DG scheme {depended} on a local continuous spacetime Galerkin predictor
to achieve high order accuracy both in space and time.
In the local continuous spacetime Galerkin predictor step, a local Cauchy problem {was} solved in each cell,
based on a weak formulation of the original partial differential equations in spacetime.
Then the high order accuracy in space and time {could} be obtained
by using the resulting spacetime representation of the numerical solution in each cell.
Our scheme {was} formulated in modal space, and the volume integral and the numerical fluxes terms at the cell interfaces in the scheme
{could} be explicitly expressed to save computational cost. This paper {provided}
the implementation details of the scheme
on two-dimensional structured meshes at third order.
The computational complexity of the ADER-DG scheme {was} compared to that of the RKDG scheme,
and extensively numerical experiments
{were} presented to show that the scheme {could} capture the viscosity solutions of the HJ equations accurately and it {was} more efficient.
{By the way, this scheme does work on unstructured grid.}

\noindent {\bf Acknowledgments.}
This work was partially supported by the Special Project on High-performance Computing
under the National Key R\&D Program (No. 2016YFB0200603), Science Challenge Project
(No. JCKY2016212A502), and the National Natural Science Foundation of China (Nos.
91330205, 91630310, 11421101).

\end{document}